\documentclass[review]{elsarticle}

\usepackage{lineno,hyperref}
\usepackage{float} 
\usepackage{arydshln}
\usepackage{color}	
\definecolor{gainsboro}{rgb}{0.86, 0.86, 0.86}
\modulolinenumbers[5]
 \usepackage{natbib}
\journal{Journal of \LaTeX\ Templates}
\usepackage{amsthm}
\newtheorem{theorem}{Theorem}
\usepackage[section]{placeins}

\usepackage{amsmath,amssymb}
\usepackage{url}
\usepackage{graphicx}
\usepackage{float} 
\usepackage{multirow}
\usepackage{color, colortbl}
\definecolor{b}{RGB}{50,153,204}
\definecolor{RawSienna}{RGB}{229, 108, 76}
\definecolor{Gray}{gray}{0.9}
\usepackage{booktabs}
\usepackage{float} 
\usepackage[super]{nth}
\newtheorem{definition}{Definition}

\begin{document}

\begin{frontmatter}

\title{Improving preference disaggregation in multicriteria decision making: incorporating time series analysis and a multi-objective approach}


\author[mymainaddress]{Betania Silva Carneiro Campello \fnref{myfootnote} \corref{mycorrespondingauthor}}
\author[secondadress]{Sarah BenAmor}
\author[thirdadress]{Leonardo Tomazeli Duarte}

\author[mymainaddress]{Jo\~ao Marcos Travassos Romano}

\address[mymainaddress]{School of Electrical and Computer Engineering (FEEC), University of Campinas (UNICAMP), Campinas SP, Brazil}
\address[secondadress]{Telfer School of Management, University of Ottawa (uOttawa), Ottawa ON, Canada}
\address[thirdadress]{School of Applied Sciences (FCA), University of Campinas (UNICAMP), Limeira SP, Brazil}
\fntext[myfootnote]{E-mail addresses: betacampello@gmail.com (B.S.C. Campello), benAmor@telfer.uottawa.ca  (S. BenAmor) leonardo.duarte@fca.unicamp.br (L.T. Duarte), romano@dmo.fee.unicamp.br (J.M.T. Romano)}

\fntext[myfootnote]{This work was supported by the National Council for Scientific and Technological Development (CNPq, Brazil), grant numbers, 168968/2018-5, 312228/2020-1.This work was also supported by the Sao Paulo Research  Foundation (FAPESP), grant 2020/09838 and grant 2020/01089-9.}

%
\cortext[mycorrespondingauthor]{betania@decom.fee.unicamp.br}

\begin{abstract}
	
	 Preference disaggregation analysis (PDA) is a widely used approach in multicriteria decision analysis that aims to extract preferential information from holistic judgments provided by decision makers. This paper presents an original methodological framework for PDA that addresses two significant challenges in this field.  Firstly, it considers the multidimensional structure of data to capture decision makers' preferences based on descriptive measures of the criteria time series, such as trend and average. This novel approach enables an understanding of decision makers' preferences in decision-making scenarios involving time series analysis, which is common in medium- to long-term impact decisions.  Secondly, the paper addresses the robustness issue commonly encountered in PDA methods by proposing a multi-objective and Monte Carlo simulation approach. This approach enables the consideration of multiple preference models and provides a mechanism to converge towards the most likely preference model. The proposed method is evaluated using real data, demonstrating its effectiveness in capturing preferences based on criteria and time series descriptive measures. The multi-objective analysis highlights the generation of multiple solutions, and, under specific conditions, reveals the possibility of achieving convergence towards a single solution that represents the decision maker's preferences.
 
\end{abstract}

\begin{keyword}
 multiple criteria decision-making \sep MCDA \sep  preference learning \sep  preference disaggregation analysis \sep PDA \sep time series analysis
\end{keyword}

\end{frontmatter}


\section{Introduction}
\label{sec:introduction}

Multiple Criteria Decision Aid (MCDA) aims to model relations between a set of alternatives, evaluated according to a set of criteria~\citep{figueira2005multiple, cinelli2020support}. A well-known MCDA approach is the preference disaggregation analysis (PDA) that infers preferential information from examples of choices provided by the decision-maker. The preference disaggregation analysis in MCDA is usually based on the additive utility function models, which have good interpretability for the decision-maker~\cite{kadzinski2017expressiveness, guo2021hybrid}. PDA is related to the field of preference analysis, which in turn relies on more typical machine learning strategies~\cite{doumpos2019preference}. However, while machine learning is data-driven, MCDA provides for greater interaction with the decision-maker~\cite{martyn2023deep}. 

Several learning approaches have been proposed in MCDA to infer the utility functions of criteria using ordinal regression techniques. Among them, UTA~\citep{jacquet1982assessing} and UTASTAR~\citep{siskos1985utastar} are two well-known methods. 


Usually, in MCDA methods, each criterion is associated with a single value that can be determined by averaging criterion performance over a period of time, using the most recent available criterion value, or relying on other static data. In this regard, the data is often represented in a matrix structure, in which the rows represent the alternatives, the columns represent the criteria, and each element represents the performance of an alternative in a certain criterion. By applying a disaggregation method to this matrix, the output yields utility functions for each criterion, indicating the decision-maker's preferences towards different criterion values.

While some decisions can be represented in the aforementioned manner, many others require considering additional information described in a multidimensional structure. For instance, certain decision-making processes require analyzing the evolution of criteria over time (time series), multiple decision-makers' preferences, or various criteria dimensions, such as social, economic, environmental, and political perspectives. Indeed, there is a growing trend in machine learning to use multidimensional data structures like \textit{tensors} or \textit{multi-arrays} to learn preferences and improve the accuracy and comprehensiveness of decision-making. For example, \cite{nilashi2021travellers} used machine learning techniques to learn user preferences on tourism websites by considering three dimensions: criteria, users, and items. Additionally, temporal analysis of user-item interactions aiming to detect temporal behavior patterns has been studied in~\cite{cho2021learning}. \cite{frolov2017tensor} conducted a review of tensor factorization techniques and their applicability in preference-based recommendations.
 
In the MCDA field, some works have emphasized the importance of incorporating additional information, such as time series data into the approaches~\cite{xu2008multi, kandakoglu2019multicriteria, martins2021multidimensional, campello2023exploiting, pinar2022choquet}. However, despite the significance of considering multidimensional structures to capture more relevant information for decision-making, only a few studies in the field of MCDA have addressed this higher-order structure~\cite{kandakoglu2019multicriteria, martins2021multidimensional, campello2023exploiting}. 

\begin{figure}[h]
	\centering
	\includegraphics[width=10cm]{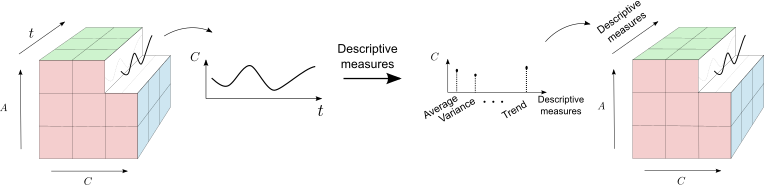}
	\caption{Tensorial structure representing the time series and descriptive measures of the criteria over time.}
	\label{fig:tensor}
\end{figure}

Given this gap, we propose a novel approach in MCDA for preference analysis that considers temporality through a tensorial structure. We introduce an extension of the UTASTAR method, which we call UTASTAR-Tensorial (or simply UTASTAR-T), to model the decision-maker's preferences related to the criteria, as well as preferences related to descriptive measures of these criteria time series (such as the criteria's time series average, trend, and seasonality). Figure~\ref{fig:tensor} illustrates the tensorial structure employed in this approach. 


In PDA, a challenge that arises when using indirect preference information is the existence of multiple compatible preference models, which result in several possible solutions that satisfy the decision-maker's preferences~\cite{kadzinski2022review}. This challenge raises a robustness issue, which has been addressed in PDA in three main streams~\cite{doumpos2019preference}. One research stream is to use simulation techniques to describe the set of feasible models, while another stream focuses on formulating robust recommendations using multiple acceptable decision models. A third stream is to develop techniques for selecting the most representative model from the set of all models compatible with the reference set information.

To address the robustness issue, we propose a novel multi-objective optimization together with Monte Carlo simulation. Using a multi-objective optimization approach, one can identify multiple compatible decision models, which provide the decision-maker with a comprehensive understanding of their preferences \cite{corrente2013robust}.

In summary, this paper's contribution is a novel methodology for analyzing decision-maker preferences by incorporating descriptive measures of criteria time series and using a multi-objective approach to generate multiple solutions. This is an original methodological contribution because we extend the UTASTAR method to a multi-array data representation and present a new approach using a multi-objective optimization method.

The paper is organized as follows: Section~\ref{sec:Study_motivation} discusses the motivation behind considering temporality in preferential learning in MCDA. Section~\ref{sec:literature_review} provides a overview of relevant works related to PDA. Section~\ref{sec:background} provides background on classical UTA and UTASTAR, as well as multi-objective problems. Section~\ref{sec:UTASTAR-T} introduces our proposed approach, the UTASTAR-T method, and the multi-objective optimization.  Section~\ref{sec:computationalresults} presents and analyzes the results. Finally, Section~\ref{sec:conclusions} concludes the paper with our final remarks.

\section{Temporal analysis in multicriteria decision making: enhancing preferential learning}
\label{sec:Study_motivation}


Time series analysis has long been used in various fields to aid decision-making by identifying patterns and extracting descriptive measures such as trends, volatility, averages, and seasonality. For example,~\cite{teichgraeber2022time} reviews studies of renewable energy that consider temporality and argue that renewable energy availability varies with time, and optimal investment solutions depending on the details of temporal data. Furthermore,~\cite{dash2021intelligent,petropoulos2022covid} show how time series analysis is useful in modeling and predicting the behavior of COVID-19; they discuss that incorporating temporal information is crucial for decision-making and policy implementations. In the economic, financial, and business fields, time series analysis is well-established. Studies such as~\cite{zhang2017signal, sung2000strategic, atella2003investment, bush2021uncertainty, henriques2011effect, bresser2014developmental, alaali2020effect} have discussed the impact of exchange rate and raw material volatility on domestic and international investment decisions. Additionally, time series analysis has been applied in a wide range of fields to support decisions including healthcare, with studies such as~\cite{rezaei2016determinants, kassakian2016clinical, mehrmolaei2018tsp, holland2020interrupted}; it has also been utilized in tourism and marketing, as evidenced by studies such as~\cite{jeffrey1988temporal, butler2001seasonality, jermsittiparsert2019behavior}, among others.

As many decisions rely on temporal analysis, it is crucial for preference learning studies to incorporate temporality. In the field of machine learning, there has been research exploring the incorporation of temporal information for learning preferences in a multi-dimensional data structure. For instance, \cite{kharfan2021data}, uses machine learning to forecast product demands by accounting for factors such as volatility, product lifecycle, and varying customer preferences. Others, like \cite{cho2021learning}, focus on temporal analysis of user-item interactions to detect behavior patterns, while \cite{zhang2021dynamic} consider time-dependent features to capture changing user interests through tensor modeling.

In PDA, the criteria involved in decision-making can be time-dependent, and analyzing such data can provide essential insights for decision-making. For instance, in investment decision-making, when ranking countries based on exchange rates and raw material prices, it may be necessary to analyze the trend, average, and volatility of these criteria \cite{campello2023exploiting}. Additionally, while determining the budget ranking of hospital departments, it can be relevant to consider the analysis of the trend and average of both patient volume and expected reimbursement criteria.

Temporal analysis is not commonly employed in PDA. However, some studies have addressed the use of time-series. For instance,~\cite{angelopoulos2019disaggregating}, which proposes time series disaggregation models to forecast electricity demand that includes an ordinal regression analysis. The goal was to create a robust additive value model consistent with historical data. Also, \cite{thesari2019decision} developed a methodology for learning the decision-maker's preferences using the time-series of the criteria. Firstly, the authors identified city departments and collected historical and statistical data on public resource distribution, creating a timeline of budget data. Secondly, the UTASTAR method was applied to determine department importance by assigning weights based on historical data, revealing long-term decision behaviors of managers.

Therefore, this study aims to model decision-maker's preferences using a tensorial structure, by considering the alternatives and criteria dimensions, and incorporating a third dimension that accounts for descriptive measures of the time series of the criteria. To the best of our knowledge, no previous studies in PDA have investigated the use of a tensor to simultaneously learn preferences related to both criteria and their time series descriptive measures.

\section{Preference disaggregation analysis and robustness issues in MCDA: an overview}
\label{sec:literature_review}

Preference disaggregation analysis involves constructing additive value function models. The PDA technique was first introduced by~\citep{jacquet1982assessing} through the development of the UTA method. Numerous extensions and variants of the UTA method have been proposed, which are commonly known as UTA-based methods. One of the variants that is considered an improvement of the UTA is the UTASTAR method~\citep{siskos1985utastar}, since it introduces two error variables for the alternatives. Since the UTA and UTASTAR methods are used for ranking problems,~\citep{zopounidis2002multicriteria} presented the UTADIS method, which is used to infer decisions related to ordered categories instead of ranking. Another variant of the UTA is the UTA-poly method~\citep{sobrie2018uta}, which assumes polynomials and splines instead of piecewise linear functions. In addition,~\cite{angilella2010non} developed a non-additive robust ordinal regression in which the utility is evaluated using the Choquet integral. Recent studies have employed UTA-based methods in conjunction with machine learning techniques, as demonstrated in studies such as~\cite{martyn2023deep, guo2021hybrid}. A comprehensive review of UTA-based methods can be found in~\cite{doumpos2022preference, kadzinski2017expressiveness}.

Several studies have applied PDA in real-world scenarios, and we cite some recent ones. Preference disaggregation based on UTASTAR was proposed by \cite{ghaderi2015understanding} and applied to a database of ranked brands to investigate the role of color in customer perception. In the marine environment, the UTASTAR method was employed by \cite{stavrou2017multicriteria} to select the optimal area for ship-to-ship cargo transfer. A SIM-UTA method was developed in \cite{manolitzas2019sim}, and the approach was applied to a Greek hospital as a case study; the results indicated that the total length of stay was the most crucial factor in the emergency department's operations. In the context of an e-commerce platform,  \cite{guo2020consumer} utilized text-mining techniques to assist product managers in identifying relevant criteria and determining their relative weight importance for facilitating multiple criteria decision-making. \cite{ehsanifar2021utastar} applied UTASTAR  to a case study involving selecting the optimal location for constructing the central warehouse of the Damghan Steel Company in Iran. \cite{gehrlein4253515active} developed a UTA-based procedure to learn the actors' preferences in the context of blockchain transactions.  In \cite{wu2023value}, the authors conducted a preference disaggregation analysis of uncertain pairwise comparisons; they applied their approach to laptop recommendations, using data collected from Amazon.com. \cite{yang2023learning} proposed a method based on the attitudinal Choquet integral within an aggregation-disaggregation paradigm to acquire individual consumer preferences on product attributes, and it was applied through real cases on \textit{TripAdvisor}.

Moreover, the robustness problem in decision models has attracted many studies that focus on both analytical techniques (also known as ROR - Robust Ordinal Regression) and simulation (also known as SOR - Stochastic Ordinal Regression) \cite{kadzinski2016scoring}. In the context of analytical techniques, post-optimization is used to identify some of the points of the polyhedron formed by the constraints of the UTA and UTASTAR models. From the solutions found, it is possible to present the average of the obtained results as the final solution \cite{jacquet1982assessing}, as well as to select the most representative decision model \cite{doumpos2007regularized}, and provide recommendations based on the various decision models obtained \cite{greco2008ordinal}. In addition, techniques have been proposed to calculate an index of stability of the obtained decision models \cite{grigoroudis2002preference} or other measures of robustness \cite{doumpos2016data}. On the other hand, SOR techniques apply Monte Carlo simulation to obtain samples of utility functions. Works such as \cite{corrente2016inducing, kadzinski2016scoring, angilella2015stochastic} provide a probability distribution of utility functions compatible with the decision-maker's preferences.

In general, numerous robustness analysis methods have been proposed. For example, \cite{grigoroudis2002preference} introduced the MUSA method which aims to infer an additive collective value function and a set of satisfaction functions to obtain value functions of these satisfaction functions. The authors proposed a measure to evaluate robustness, the normalized standard deviation of the results obtained with the model. In addiction, \cite{doumpos2016data} focused on robustness metrics in a data-driven perspective, i.e., in terms of the variations in the data instances used in classification problems.~\cite{kadzinski2022review} conducted a review and experimental comparison to develop a unique recommendation in a preference disaggregation setting.

\cite{greco2008ordinal} proposed a generalization of the UTA method, called $UTA^{GMS}$, which provides a set of value functions compatible with the preferences instead of presenting a single value function in ranking problems. \cite{greco2010multiple} proposed the $UTA^{GMS}$ method for the case of sorting problems. \cite{tervonen2009stochastic} proposed the SMAA-TRI method for handling imprecise preference information in sorting problems. In addition, \cite{corrente2016inducing} introduced a series of steps to provide a probability distribution of the value functions compatible with the decision-maker's preferences, using the HAR method~\cite{smith1984efficient} to obtain these probability functions.

In \cite{kadzinski2016scoring}, scoring procedures were proposed to transform robust results into a single recommendation, using both robust and stochastic ordinal regression. In \cite{kadzinski2013stochastic}, ROR and SMAA were combined for MCDA sorting problems, and the approach was applied to classify 27 countries into four political regimes: full democracies, flawed democracies, hybrid regimes, and authoritarian regimes. In \cite{kadzinski2017expressiveness}, the authors evaluate the robustness of preference disaggregation through additive value functions in an experimental way. They presented measures used to quantify robustness and proposed a new supervised algorithm for discretization in pairwise comparisons. An illustrative example was presented to validate the algorithm.

In brief, both ROR and SOR methods provide decision-makers with multiple preference models that are compatible with the given decision scenario. As these approaches involve a constructive learning process, decision-makers typically participate in the preference learning process. Therefore, these models are useful in helping decision-makers understand their preferences and in providing recommendations for decision making \cite{corrente2013robust}. 

Finally, we have not found  any study proposing to generate multiple solutions using a multi-objective approach.

\section{Background on classical UTA and UTASTAR, and multi-objective problems}
\label{sec:background}

In Subsection~\ref{subsec:background_UTASTAR}, an introduction to classical UTA and UTASTAR is provided, while Subsection~\ref{subsec:background_multiobjectiveproblems} offers an overview of multi-objective problems.

\subsection{Classical UTA and UTASTAR: notation and linear program}
\label{subsec:background_UTASTAR}


This section briefly reviews two classical methods based on the additive value function model, the UTA~\cite{jacquet1982assessing} and UTASTAR~\cite{siskos1985utastar} methods. These methods are employed to learn the decision-maker's preferences from known examples of partial decisions. The UTA methods family infer that such a preference can be represented using an additive value function,  named \textit{marginal value} or \textit{utility function},  here represented by $u_j$, $j = 1,\ldots, n$. 

Moreover, following the general modeling methodology of MCDA problems, let us consider a set of reference alternatives (actions), $A_R = \{a_1,\ldots, a_m\}$, and a
consistent family of criteria, $C = \{c_1, \ldots, c_n\}$. Each criterion is a non-decreasing real valued function defined on A, as~\citep{zopounidis2010handbook}:
\begin{equation}
	c_j: A_R \rightarrow [c_{j*}, c_{j}^*] \subset \mathbb{R}/a \rightarrow c(a) \in \mathbb{R} 
\end{equation}

\noindent where $[c_{j*}, c_{j}^*]$ is the criterion evaluation scale, and $c_{j*}$ and $c_{j}^*$ represent the worst and best level of the $j$-th criterion. When the evaluation scales of some $c_j \in C$ are continuous (or when they contain a considerable  number of grades), the range $[c_{j*}, c_{j}^*]$ can be split into ($\alpha_j - 1$) equal intervals; otherwise,  $\alpha_j$ is similar to the number of grades. The criterion evaluation scale of each criterion $j$ is then represented by $[c_{j}^1, c_{j}^2, \cdots, c_{j}^{\alpha_j}]$, where $c_{j}^1 = c_{j*}$ and $c_{j}^{\alpha_j} = c_{j}^*$.


The global value function of  alternative $i$ is represented by $u[\textbf{c}(a_i)]$, $i =1,\ldots, m$, in which the following properties hold:
\begin{equation} 
	\begin{cases} 
		&	u[\textbf{c}(a_i)] > u[\textbf{c}(a_{i+1})] \iff a_i \succ a_{i+1} \ \ \mbox{(preference)}  \\
		&	u[\textbf{c}(a_i)] = u[\textbf{c}(a_{i+1})] \iff a_i \sim a_{i+1} \ \ \mbox{(indifference)}
	\end{cases}\label{eq:comparison}
\end{equation}

Linear programming is used to estimate the decision-maker's preference model. To this end, firstly, the set of reference alternatives, $A_R = \{a_1,\ldots, a_m\}$, is ordered in such a way that $a_1$ is the best and $a_m$ the worst alternative according to the decision-maker's ranking. Then, each pairwise comparison showed in~(\ref{eq:comparison}) is translated into a constraint in the UTA model as following: for each pair of alternatives $(a_i, a_{i+1}) \in A_R$, we have $u[\textbf{c}(a_i)] - u[\textbf{c}(a_{i+1})] > 0$ iff $a_i \succ a_{i+1}$, and $u[\textbf{c}(a_i)] - u[\textbf{c}(a_{i+1})] = 0$ iff  $a_i \sim a_{i+1}$. However, note that these constraints can be infeasible; therefore, error variables, $\sigma(a)$, are introduced, and we have:
\begin{equation} 
	\Delta(a_i, a_{i+1}) = u[\textbf{c}(a_i)] - u[\textbf{c}(a_{i+1})] + \sigma(a).\label{eq:prefindif}
\end{equation}

\noindent The UTA constraints become:
\begin{equation} 
	\begin{cases}
		&  \Delta(a_i, a_{i+1}) \geq \delta \ \mbox{if} \ a_i \succ a_{i+1} \ \forall i \\
		&  \Delta(a_i, a_{i+1}) = 0 \ \mbox{if} \ a_i \sim a_{i+1}  \ \forall i 
	\end{cases}\label{eq:constraint},
\end{equation}
\noindent where $\delta$ is a small positive number so as to discriminate significantly two
successive equivalence classes of $R$. If all the error variables are zero, the decision-maker's preference order is feasible; therefore, to achieve this preference order, the objective function of UTA aims to minimize the error variables. The linear program is given by:

\begin{align}
	min \ z =  & \ \sum_{a \in A_r} \sigma(a) \\
	\mbox{subject to} \nonumber \\
	&  \Delta(a_i, a_{i+1}) \geq \delta \ \mbox{if} \ a_i \succ a_{i+1}, \ \forall i  \label{eq:uta-DELTA}  \\
	&  \Delta(a_i, a_{i+1}) = 0 \ \mbox{if} \ a_i \sim a_{i+1}, \ \forall i   \label{eq:deltaindiference} \\
	& u_j(c_j^{\ell + 1}) - u_j(c_j^\ell) \geq 0, \ \  j = 1,\ldots, n, \ell = 1, \ldots, \alpha_j  \label{eq:monotonicity} \\
	&  \sum_{j=1}^{n} u_j(c_{j}^*) = 1 \label{eq:equaltoone} \\
	& u_j(c_{j*}) = 0, u_j(c_j^\ell) \geq 0, \sigma(a) \geq 0, \ \forall a \in A_R, \ \forall \ j, \ell \label{eq:uta-BOUNDS}
\end{align}

\noindent Finally, Constraint~(\ref{eq:monotonicity}) guarantees the monotonicity of the utility functions, and Constraint~(\ref{eq:equaltoone}), is related to the normalization of these functions.

UTASTAR is considered an improved version of UTA~\citep{siskos2005uta, zopounidis2010handbook}. Two important modifications were proposed. First, a double error function, $\sigma^+(a_i)$ and  $\sigma^-(a_i)$, is introduced, which better minimizes the dispersion of points around the monotone utility function curve. Then, Equation~(\ref{eq:prefindif}) is rearranged as
\begin{equation} 
	\Delta(a_i, a_{i+1}) = \left(u[\textbf{c}(a_i)] - \sigma^+(a_i) +  \sigma^-(a_i) \right) - \left(u[\textbf{c}(a_{i+1})] - \sigma^+(a_{i+1}) +  \sigma^-(a_{i+1}) \right) \label{eq:prefindifutastar}
\end{equation}

\noindent Second, Constraint~(\ref{eq:monotonicity}) is replaced by the new transformation variables, $w_{j \ell}$, which reduces the size of the LP and represent the successive steps of the utility functions, in the form
\begin{equation}
	w_{j \ell} = u[\textbf{c}^{\ell + 1}(a_i)] - u[\textbf{c}^{\ell}(a_i)] \geq 0,
\end{equation}

\noindent for simplicity's sake, we represent it as
\begin{equation}
	w_{j \ell} = u(\textbf{c}^{\ell + 1}) - u(\textbf{c}^{\ell}) \geq 0.
\end{equation}

\noindent Thus, the linear program of the UTASTAR is given by:	
\begin{align}
	min \ z =  & \ \sum_{i=1}^{m} [\sigma^+(a_i) +  \sigma^-(a_i) ] \\
	\mbox{subject to} \nonumber \\
	&  \Delta(a_i, a_{i+1}) \geq \delta \ \mbox{if} \ a_i \succ a_{i+1} \ \forall i \\
	&  \Delta(a_i, a_{i+1}) = 0 \ \mbox{if} \ a_i \sim a_{i+1}  \ \forall i \\
	&  \sum_{j=1}^{n} \sum_{\ell = 1}^{\alpha_j - 1} w_{j \ell} = 1,  \label{eq:wsumone} \\
	& w_{j \ell} \geq 0, \sigma^+(a_i) \geq 0,  \sigma^-(a_i) \geq 0, \ \forall \ j, \ell \  \mbox{and} \  i 
\end{align}

\subsubsection{Robustness issue in UTA and UTASTAR}
\label{subsubsec:RobustnessprobleminUTAandUTASTAR}

The UTA and UTASTAR methods are very useful for preference analysis, however they present a robustness issue in the decision model. In \cite{jacquet1982assessing}, the authors emphasize that there are different feasible solutions (i.e., different utility functions) that are compatible with the decision-maker's preferences. If the possible solutions are very different from each other, the model will not be representative and will provide different recommendations when applied to a new set of alternatives \cite{doumpos2019preference}. Thus, when the optimal solution $z^*$ is equal to zero, different solutions can model the decision-maker's preferences. Furthermore, when the optimal solution is different from zero, it is possible to increment $z^*$ in the form of $z^* + \ (z^*)$, where $\gamma$ is a very small value, and still obtain different solutions. Thus, \cite{jacquet1982assessing} proposed a post-optimization analysis by solving the following linear program. For $j =1, \ldots, n$:
\begin{equation}\label{eq:minmax}
	\begin{cases}
		& \mbox{min} \ \ u_j(c_j^*) \\		
		& \mbox{subject to}: \\
		& \mbox{Constraints} \ (\ref{eq:uta-DELTA})-(\ref{eq:uta-BOUNDS}) \\
		& z \leq z^* + \gamma (z^*)
	\end{cases}\ \ \ \mbox{and} \ \ \ \ 
	\begin{cases}
		& \mbox{max} \ \ u_j(c_j^*) \\		
		& \mbox{subject to}: \\
		& \mbox{Constraints} \ (\ref{eq:uta-DELTA})-(\ref{eq:uta-BOUNDS}) \\
		& z \leq z^* + \gamma (z^*)
	\end{cases} 
\end{equation}

The final solution is the average of the solutions obtained.

\subsection{Multi-objective problems}
\label{subsec:background_multiobjectiveproblems}

Optimization problems consist of identifying, within a set of solutions, the solution that is the best, given the objective function and the constraints of the problem. Multi-objective optimization is a mathematical technique used to find optimal solutions in problems where there are multiple conflicting objectives to be optimized~\cite{coello2007evolutionary, miettinen2012nonlinear,campello2020multiobjective}. Due to the conflicting nature of the objective functions, it is not possible to find a single solution that is optimal for all objectives simultaneously, which leads to multiple optimal solutions. The set of potential solutions is the \textit{decision space}, represented here by $X$, and its elements $x$ are the feasible solutions.  Each feasible solution can be mapped into its corresponding objective vector $y = f(x)$ and $Y = \{y \in  \mathbb{R}^n: y = f(x) \ \mbox{for some} \ x \in X\}$ which is the set of feasible outcomes in the objective space. 

The multi-objective optimization problem is formulated as:
\begin{align}
	&\text{min  }  f(x) \ = \{f_1(x), f_2(x), \cdots, f_n(x)\}  \label{eq:modelomultiobjetivo} \\
	&\text{s.t.} \ \  x \in X \label{eq:modelomultiobjetivo2},
\end{align}
\noindent where there are $n > 2$ objective functions, $f_j: \mathbb{R}^n \rightarrow \mathbb{R}$, and   
the optimal solutions of~(\ref{eq:modelomultiobjetivo})-(\ref{eq:modelomultiobjetivo2}) are called \textit{efficient solutions} or \textit{Pareto optimal solutions}. The optimality is a crucial concept in optimization problems, and it is defined as follows.
\begin{definition}
	A solution $x^* \in X$ is said to be an weakly efficient solution if there exists no feasible solution $x \in X$ such that $f_j(x) < f_j(x^*)$, $\forall j = 1, \ldots, n$
\end{definition}
\begin{definition}
	A solution $x^* \in X$ is said to be an efficient solution if there exists no feasible solution $x \in X$ such that $f_j(x) \leq f_j(x^*)$, $\forall j = 1, \ldots, n$, and $f_{j'}(x) < f_{j'}(x^*)$ for some $j' \in \{1, \ldots, n\}$.
\end{definition}

There are several methods to solve~(\ref{eq:modelomultiobjetivo})-(\ref{eq:modelomultiobjetivo2}), such as classical techniques, outlined in~\cite{miettinen2012nonlinear}, and evolutionary algorithms, as described in~\cite{coello2007evolutionary}. In this study, we use the classical weighting method~\cite{gass1955computational, miettinen2012nonlinear} due to its simple concept, ease of implementation, and its ability to guarantee, at least, weakly Pareto optimal solutions, as we will later discuss.

The weighting method consists of assigning weights $\mu_j$, with $\sum_{j =1}^{n} \mu_j = 1$ and $\mu_j \in [0,1]$, to the $n$ objective functions, transforming the problem into a weighted single-objective problem. Thus, to apply this method, we rewrite the multi-objective optimization~(\ref{eq:modelomultiobjetivo})-(\ref{eq:modelomultiobjetivo2}) in the following form:
\begin{align}
	\text{min  }   & f(x) = \sum_{j = 1}^{n} \mu_j f_j(x) \label{eq:somaponderada} \\
	\text{s.t}  
	& \ \ \ x \in X.\label{eq:somaponderada2}
\end{align}

\noindent Therefore, by varying the values of $\mu_j$, an optimal solution on the Pareto curve can be obtained.

By setting all $\mu_j$ to zero, except for a single $j$ (for example, $\mu_1 = 1$, and $\mu_j = 0$, for $j = 2, \cdots, n$), Equations~(\ref{eq:somaponderada})-(\ref{eq:somaponderada2}) become:
\begin{align}
	\text{min  }   & f(x) =  1 \times f_j(x), \ \ j =1, \cdots, n \label{eq:modelomultiobjetivoFO3} \\
	\text{s.t}  
	& \ \ \ x \in X, \label{eq:modelomultiobjetivoRESt3}
\end{align}  
which is equivalent to minimizing each objective function $f_j$ individually. By solving~(\ref{eq:modelomultiobjetivoRESt3}) for all $j$, one obtains the lower bound of each function, while respecting the Constraints~(\ref{eq:modelomultiobjetivoRESt3}).

Theoretical results are presented concerning the weighting method solutions. As proofs of the theorems can be found in~\cite{miettinen2012nonlinear}.

\begin{theorem}\label{theorem:weaklyParetooptimal}
	The solution of the weighting problem (\ref{eq:somaponderada})-(\ref{eq:somaponderada2}) is weakly Pareto optimal.
\end{theorem}

\begin{theorem}\label{theorem:wParetooptimal}
	The solution of the weighting problem (\ref{eq:somaponderada})-(\ref{eq:somaponderada2}) is Pareto optimal if the weighting coefficients are all positives: $\mu_j > 0, \ \ \forall j=1,\ldots, n$.
\end{theorem}

\noindent Theorems~\ref{theorem:weaklyParetooptimal} and~\ref{theorem:wParetooptimal} guarantee that the solution cannot be improved for one objective without worsening it for another objective. Therefore, the weighting method is a simple but efficient technique. 


%

\section{UTASTAR-T with a multi-objective approach}
\label{sec:UTASTAR-T}

This section presents the UTASTAR-T method, an extension of UTASTAR for a situation in which the input data is represented by a tensor $\mathcal{P} \in \mathbb{R}^{m \times n \times T}$, with $m$ alternatives, $n$ criteria, and $T$ time series samples. The steps of the methodology consist of applying a function to the tensor $\mathcal{P}$ to calculate descriptive measures of the time series, from which the tensor $\mathcal{S} \in \mathbb{R}^{m \times n \times h}$ is obtained, with $h$ descriptive measures. Then, the UTASTAR-T method is applied to the tensor $\mathcal{S}$, from which the value functions of the $n$ criteria on the $h$ descriptive measures are obtained.

Thus, in this approach, instead of estimating the value function of each criterion, $u_j(c_j)$, as in UTASTAR, we infer the value functions of the criteria for each descriptive measure $k$, denoted by $u_{jk}(c_{jk})$, $k = 1,\ldots, h$. Therefore, we have $n$ criteria and $h$ descriptive measures that constitute $n \times h$ value functions. The global value functions of alternative $i$ is represented by $u[{c}(a_i)]$, $i =1,\ldots, m$, obtained as follows:
\begin{equation}
	u[\textbf{c}(a_i)] = \sum_{k=1}^{h} \sum_{j=1}^{n} u_{kj}(c_{kj}), i =1,\ldots, m. \label{eq:globalutensor}
\end{equation}

The global value of an alternative is defined as a linear combination of the utilities for each criterion and for the different defined descriptive measures.

The value scale of each criterion $j$ for each descriptive measure $k$ is denoted by $[c_{jk}^1, c_{jk}^2$, $\cdots, c_{jk}^{\alpha_{kj}}]$, where $c_{jk}^1$ is the worst value of criterion $j$ for descriptive measure $k$ (also represented by $c_{jk*}$) and $c_{jk}^{\alpha_{kj}}$ is the best value (also represented by $c_{jk}^*$). Properties~(\ref{eq:comparison}) and Equation~(\ref{eq:prefindifutastar}) are similar to those of UTASTAR.

There are two other important modifications concerning the variables $w_{j \ell}$. The first one is that these variables are replaced by $w_{kj \ell}$ in order to ensure monotonicity conditions for all value functions in each descriptive measure, as follows:
\begin{equation}\label{eq:wemfuncaodeu}
	w_{kj \ell} = u_{kj}(c^{\ell+1}) -  u_{kj}(c^{\ell}) \geq 0, \ \forall j, \ell, k.
\end{equation}

\noindent Moreover, in classical UTASTAR, constraints~(\ref{eq:wsumone}) normalize the weights of the criteria, which should sum up to 1. In UTASTAR-T, $w_{kj \ell}$ are normalized in such a way that they sum up to 1 for each descriptive measure: $\sum_{j=1}^{n} \sum_{\ell = 1}^{\alpha_{kj} - 1} w_{kj \ell} = 1$, for each $k$.

Taking these modifications into account, the proposed UTASTAR-T algorithm can be summarized as follows.

\noindent \textbf{Step 1:}

\noindent Given the tensor $\mathcal{P} \in \mathbb{R}^{m \times n \times T}$, transform it into the  descriptive measures space:
\begin{eqnarray}
	f: \mathcal{P}  \in \mathbb{R}^{m \times n \times T} \rightarrow \mathcal{S}  \in \mathbb{R}^{m \times n \times h},
\end{eqnarray}

\noindent where $f$ is the function that maps from the domain $\mathcal{P}$ to $\mathcal{S}$, the $s_{ijk}$ are the elements of $\mathcal{S}$, and $h$ is the number of descriptive measures. The choice of descriptive measures should be made according to the decision-maker's preferences. In this study, we consider the average ($k=1$), calculated according to Equation~(\ref{average}), and the slope coefficient ($k=2$), calculated using Equation~(\ref{sc}).
\begin{eqnarray}\label{average}
	s_{ij1} = \hat{p}_{ij} = \frac{1}{T} \sum_{t=1}^{T} p_{ijt}, \ \ \forall \ i, j, 
\end{eqnarray}
where $\hat{p}_{ij}$ represents the average of the time series in alternative $i$ and criterion $j$. Trend ($k=2$) is calculated through the angular coefficient ($s_{ij2} = \hat{\beta}$) obtained by linear regression of the form: $s_{ij2} = \hat{\alpha} + \hat{\beta} t$, $\forall \ t = 1, \cdots, T$, where $\hat{\beta}$ is obtained as
\begin{eqnarray}\label{sc}
	\hat{\beta} = \frac{\sum_{1}^{T}(p_{ijt} - \hat{p}_{ij})(t_{ijt} - 	\hat{t}_{ij})}{\sum_{1}^{T}(p_{ijt} - \hat{p}_{ij})^2} , \ \ \forall \ i, j,
\end{eqnarray}

\noindent \textbf{Step 2: }

\noindent Calculate the global value of the alternatives $u[\textbf{c}(a_i)]$, $i = 1, \cdots, m$, in terms of marginal values $u_{kj}(c_{jk})$ according to Equation~(\ref{eq:globalutensor}), and then in terms of the variables $w_{kj \ell}$, as follows:
\begin{equation} 
	\begin{cases} 
		& u_{kj}(c_{kj}^1) = 0, j = 1, \ldots, n \\
		& u_{kj}(c_{kj}^\ell) = \sum_{\ell = 1}^{\alpha_{kj} -1} w_{kj \ell}, j = 1, \ldots, n \ \ \mbox{e} \ \ k =  1, \ldots, h \\
	\end{cases}	\label{eq:utilitequaw}
\end{equation}

\noindent \textbf{Step 3:}

\noindent Perform pairwise comparison of consecutive alternatives in the ranking and add the error function as in Equation (\ref{eq:prefindifutastar}).

\vspace{0.2cm}

\noindent \textbf{Step 4:}

\noindent Solve the linear program:

\begin{align}
	\min \ z =  & \ \sum_{i=1}^{m} [\sigma^+(a_i) +  \sigma^-(a_i) ] \label{eq:utastartFO} \\
	\mbox{subject to} \nonumber \\
	&  \Delta(a_i, a_{i+1}) \geq \delta \ \mbox{iff} \ a_i \succ a_{i+1} \ \forall i \label{eq:delta1} \\
	&  \Delta(a_i, a_{i+1}) = 0 \ \mbox{iff} \ a_i \sim a_{i+1}  \ \forall i \\
	&  \sum_{j=1}^{n} \sum_{\ell = 1}^{\alpha_{kj} - 1} w_{kj \ell} = 1, \ \forall \ k \\
	& w_{kj \ell} \geq 0, \sigma^+(a_i) \geq 0,  \sigma^-(a_i) \geq 0, \ \forall \ j, \ell \  \mbox{and} \  i \label{eq:utastartultimaequacao}
\end{align}
\noindent After applying UTASTAR-T, we obtain the values of each description measure for each criterion in the form of weights $w_{kj\ell}$.

\subsection{Multi-objective UTASTAR-T linear program}
\label{subsec:multiobjetctiveutastart}

In Subsection~\ref{subsubsec:RobustnessprobleminUTAandUTASTAR}, we presented the proposal of the study \cite{jacquet1982assessing} for post-optimization analysis, which consists of maximizing and minimizing each value function $u_j(c_j^*)$, $\forall \ j$, according to Equations~(\ref{eq:minmax}), and calculating the average of the obtained solutions. This process of maximizing and minimizing each value function is equivalent to obtaining the upper and lower bounds of each value function $j$, similar to what was done in Equations~(\ref{eq:modelomultiobjetivoFO3})-(\ref{eq:modelomultiobjetivoRESt3}). Therefore, instead of obtaining only the bounds of the value functions and calculating the average, we propose to apply the multi-objective weighted sum method to analyze multiple optimal solutions that the value functions can assume.

To test the existence of multiple optimal solutions, we propose the following multi-objective optimization:
\begin{align}\label{eq:multiobjectivemath}
	max \ F(x) =  & (u_{11}(c_{11}), \ldots, u_{1n}(c_{1n}),  \ldots, u_{21}(c_{21}), \ldots, u_{2n}(c_{2n}), \nonumber \\ 
	&\ldots, u_{h1}(c_{h1}), \ldots, u_{hn}(c_{hn}))^T
\end{align}

\noindent in the polyhedron of Constraints (\ref{eq:delta1})-(\ref{eq:utastartultimaequacao}) bounded by the new constraint:
\begin{align}\label{eq:novaconstraint}
	\sum_{i=1}^{m} [\sigma^+(a_i) +  \sigma^-(a_i) ] \leq z^* +\   \epsilon,
\end{align}

\noindent where z* is the optimal value obtained by solving (\ref{eq:utastartFO})-(\ref{eq:utastartultimaequacao}) in \textbf{Step 4} of the UTASTAR-T algorithm and $\epsilon$ is a very small positive number.

To solve the multi-objective optimization, we used simulation and the weighted sum method, assigning a weight $\mu_j$ to each $u_{kj}(c_{kj})$ in the multi-objective function (\ref{eq:multiobjectivemath}), and repeating the following steps 1,000 times:

\noindent \textbf{Step MO-1}: Generate random samples for each $\mu_j$, such as: $\mu_j \sim U[0, 1]$, $\forall j$;

\noindent \textbf{Step MO-2}: Solve the following Multi-objective optimization with the weighted sum method:
\begin{align}
	max \ F(x) =  &  \ \mu_1 \times \left(\sum_{k=1}^{h}u_{kj}(c_{k1})\right)  +  \mu_2 \times \left(\sum_{k=1}^{h}u_{kj}(c_{k2})\right) \nonumber \\
	& + \ldots + \mu_n \times \left(\sum_{k=1}^{h}u_{kn}(c_{kn})\right),
\end{align}

\noindent subject to Constraints (\ref{eq:delta1})-(\ref{eq:utastartultimaequacao}) and (\ref{eq:novaconstraint});

\noindent \textbf{Step MO-3}: Save the different solutions of vector $\textbf{w}_{k}$, $\forall k$, and the number of occurrences of each vector.

Once the solutions are obtained, it is possible to analyze how many times each solution $\textbf{w}_{k}$ was obtained in the simulation given the different values of $\mu_j$. It is also possible to calculate the weighted average of the solution vector values $\textbf{w}_{k}$. In addition, in MCDA disaggregation approaches, some authors, after obtaining the results, consider a possible interaction between the analyst and the decision-maker, in order to present, analyze, and discuss the obtained solutions~\cite{greco2010multiple, martyn2023deep, doumpos2022preference}. In this case, if the decision-maker provides preferences between the criteria, it is possible to apply an order of $\mu_j$ in the simulation. Therefore, the relevance of the criteria is considered according to the decision-maker's preference: $c_a \succ c_b \succ \cdots \succ c_y$, and the \textbf{Step MO-1} of the simulation becomes:
\vspace{0.2cm}

\noindent \textbf{Step MO-1 - With criteria ranking:} Randomly sample weights vector $\mathbf{\mu}$: $\mathbf{\mu} \sim U[0, 1]$ and order the vector as follows: $\mu_a > \mu_b > \cdots > \mu_y$, given that $c_a \succ c_b \succ \cdots \succ c_y$, obtaining $\mathbf{\mu}_{ordered} = [\mu_a, \mu_b, \cdots, \mu_y]$.

After applying \textbf{Step MO-1 - With criteria ranking}, we apply \textbf{Step MO-2} and \textbf{Step MO-3}, obtaining the solution vectors $\textbf{w}_{k}$, $\forall k$, and the number of occurrences of each $\textbf{w}_{k}$.

\section{Computational experiments using real data}
\label{sec:computationalresults}

In this section, we use the same data from the studies~\cite{banamar2018extension, campello2023exploiting}, as shown in Table~\ref{tab:dados}.We assume that the decision-maker's ranking of alternatives is the one shown in the table: Malaysia (MY) $\succ$ Russia (RU) $\succ$ Turkey (TR) $\succ$ Brazil (BR) $\succ$ China (CN) $\succ$ India (IN) $\succ$ Indonesia (ID) $\succ$ Mexico (MX) $\succ$ Philippines (PH) $\succ$ South Africa (ZA). We also assume that the decision-maker considered two relevant descriptive measures to make this decision, the mean and the slope coefficient.

\begin{table}[htbp]
	\centering
	\caption{Economic and social indicators of ten emerging countries, data obtained from~\cite{banamar2018extension}.}
	\scalebox{0.6}{\begin{tabular}{|c|r|r|r|r|r|r|r|r|r|r|r|r|r|r|r|r|r|r|}
			\hline
			Weights \textbf{w} & \multicolumn{ 6}{c|}{0.333} & \multicolumn{ 6}{c|}{0.333} & \multicolumn{ 6}{c|}{0.333} \\ \hline
			Max/Min & \multicolumn{ 6}{c|}{Max} & \multicolumn{ 6}{c|}{Max} & \multicolumn{ 6}{c|}{Max} \\ \hline
			Criteria & \multicolumn{ 6}{c|}{Life expectancy at birth ($c_1$)} & \multicolumn{ 6}{c|}{Education ($c_2$)} & \multicolumn{ 6}{c|}{Gross national income per capita ($c_3$)} \\ \hline
			\rowcolor{Gray}
			Anos  & 1990 & 1995 & 2000 & 2005 & 2010 & 2015 & 1990 & 1995 & 2000 & 2005 & 2010 & 2015 & 1990 & 1995 & 2000 & 2005 & 2010 & 2015 \\ \hline
			BR & 65.3 & 67.6 & 70.1 & 71.9 & 73.3 & 74.8 & 8.00 & 8.95 & 9.95 & 10.15 & 11.05 & 11.60 & 10065 & 10959 & 11161 & 12032 & 14420 & 15062 \\ \hline
			CN & 69.0 & 69.9 & 71.7 & 73.7 & 75.0 & 76.0 & 6.80 & 7.25 & 7.85 & 8.75 & 9.85 & 10.30 & 1520 & 2508 & 3632 & 5632 & 9387 & 13347 \\ \hline
			IN & 57.9 & 60.4 & 62.6 & 64.5 & 66.5 & 68.4 & 5.35 & 5.90 & 6.45 & 7.35 & 8.25 & 8.55 & 1754 & 2046 & 2522 & 3239 & 4499 & 5814 \\ \hline
			ID & 63.3 & 65.0 & 66.3 & 67.2 & 68.1 & 69.1 & 6.75 & 7.20 & 8.70 & 9.30 & 9.95 & 10.30 & 4337 & 5930 & 5308 & 6547 & 8267 & 10130 \\ \hline
			MY & 70.7 & 71.8 & 72.8 & 73.6 & 74.1 & 74.8 & 8.10 & 8.90 & 10.25 & 10.15 & 11.35 & 11.35 & 9772 & 13439 & 14500 & 17157 & 19725 & 23712 \\ \hline
			MX & 70.8 & 72.8 & 74.4 & 75.3 & 76.1 & 77.0 & 8.05 & 8.55 & 9.15 & 9.85 & 10.50 & 10.80 & 12074 & 12028 & 14388 & 14693 & 15395 & 16249 \\ \hline
			PH & 65.3 & 66.1 & 66.7 & 67.2 & 67.7 & 68.3 & 8.70 & 8.95 & 9.50 & 9.75 & 9.75 & 10.20 & 3962 & 4111 & 4994 & 6058 & 7478 & 8232 \\ \hline
			RU & 68.0 & 66.0 & 65.1 & 65.8 & 68.6 & 70.3 & 10.95 & 10.85 & 11.85 & 12.60 & 13.10 & 13.35 & 19461 & 12011 & 12933 & 17797 & 21075 & 22094 \\ \hline
			ZA & 62.1 & 61.4 & 55.9 & 51.6 & 54.5 & 57.9 & 8.95 & 10.65 & 11.00 & 11.15 & 11.55 & 11.75 & 9987 & 9566 & 9719 & 10935 & 11833 & 12110 \\ \hline
			TR & 64.3 & 67.0 & 70.0 & 72.5 & 74.2 & 75.6 & 6.70 & 7.20 & 8.30 & 8.95 & 10.55 & 11.05 & 10494 & 11317 & 12807 & 14987 & 16506 & 18976 \\ \hline
	\end{tabular}}
	\label{tab:dados} 
\end{table}

\begin{table}[h]
	\centering
	\caption{Results of the variables $w_{kj\ell}$ obtained after applying UTASTAR-T.}
	\begin{tabular}{ccl}
		\toprule
		Descriptive measure & $w_{kj\ell}$ & Value \\ \hline
		\multirow{3}{*}{Average} & $ w_{113}$ & 0,10\\
		&$w_{137}$  & 0,05\\
		&$w_{138}$  & 0,85\\
		\midrule 
		\multirow{4}{*}{Slope coefficient} & $w_{212}$  & 0,60\\
		& $w_{217}$  & 0,20\\
		& $w_{222}$  & 0,05\\
		& $w_{234}$  & 0,15
		\\\hline
	\end{tabular}\label{tab:variablesresults}
\end{table}
\begin{figure}
	\centering
	\includegraphics[width=6cm]{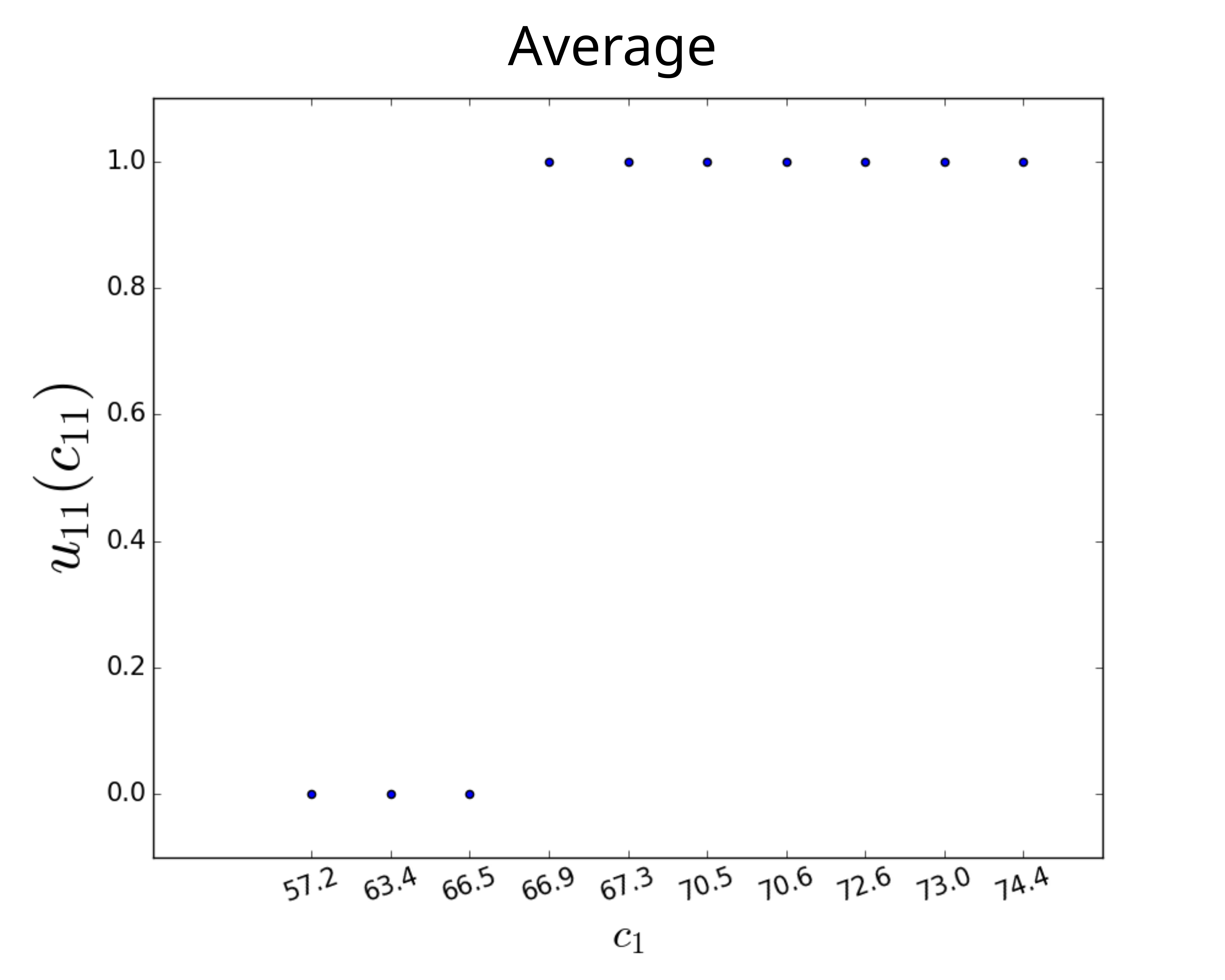}
	\includegraphics[width=6cm]{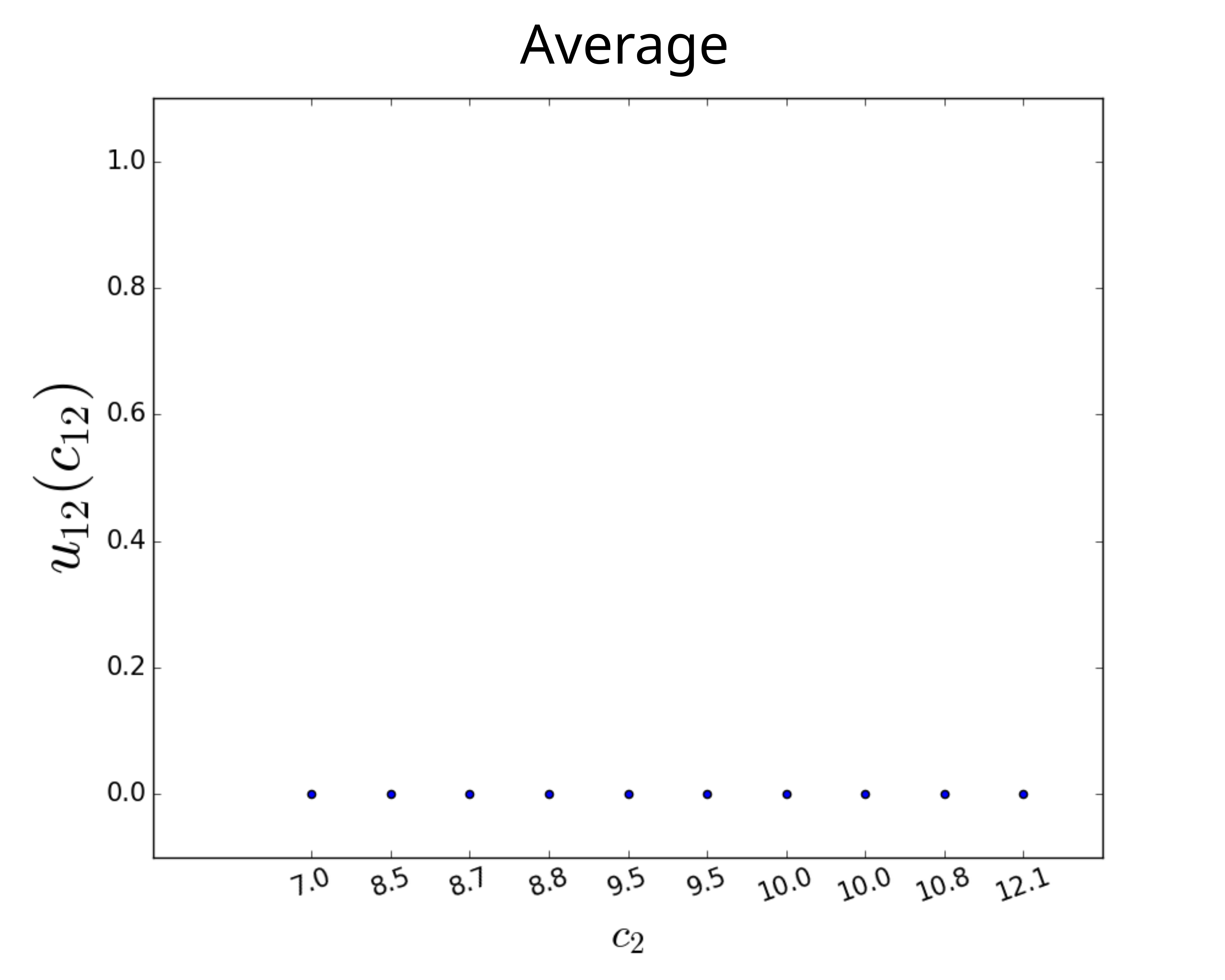}
	\includegraphics[width=6cm]{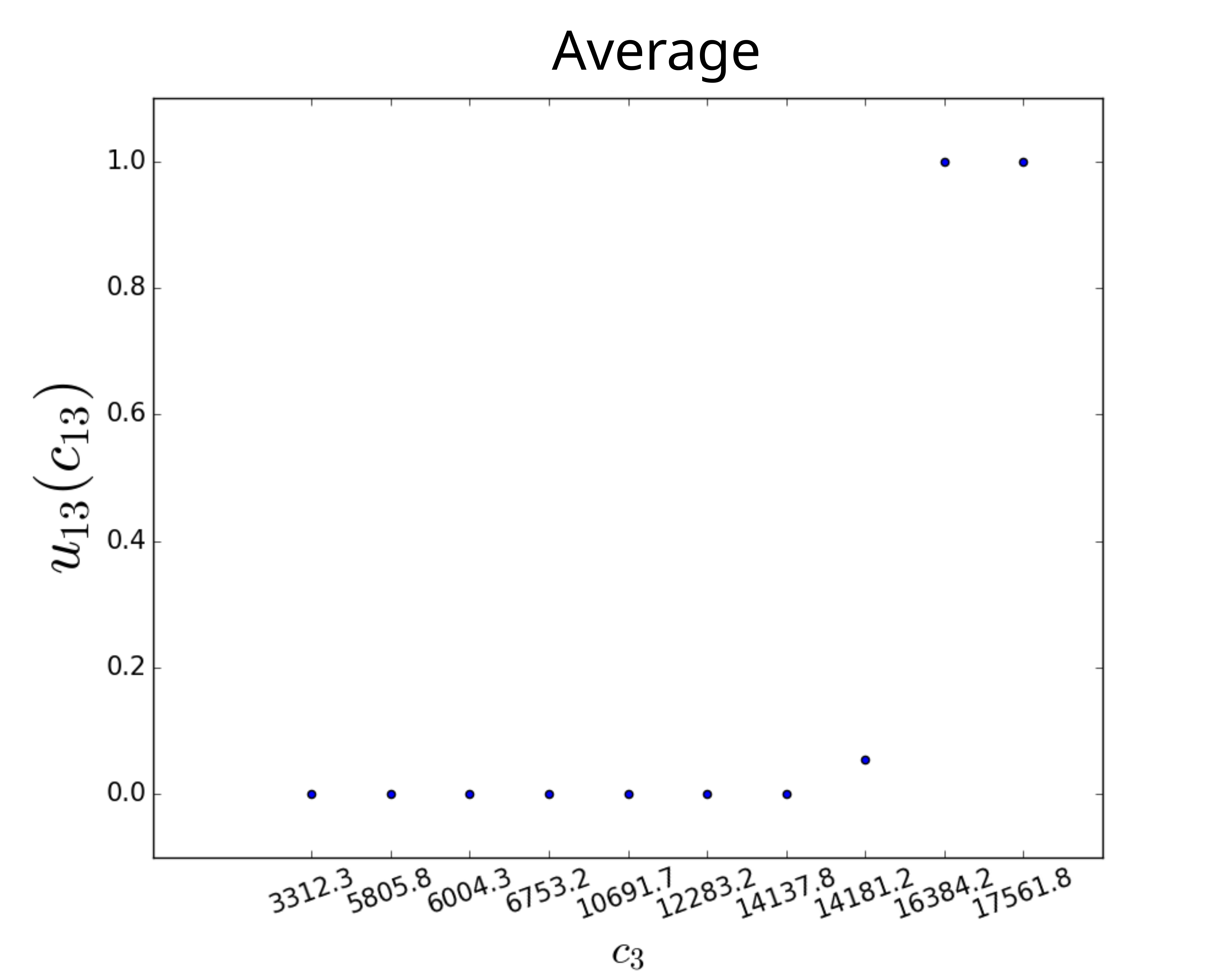}
	\caption{Normalized value functions of descriptive measure 1: average.}\label{fig:k1utilitygraph}
\end{figure}	
\begin{figure}[h]
	\centering
	\includegraphics[width=6cm]{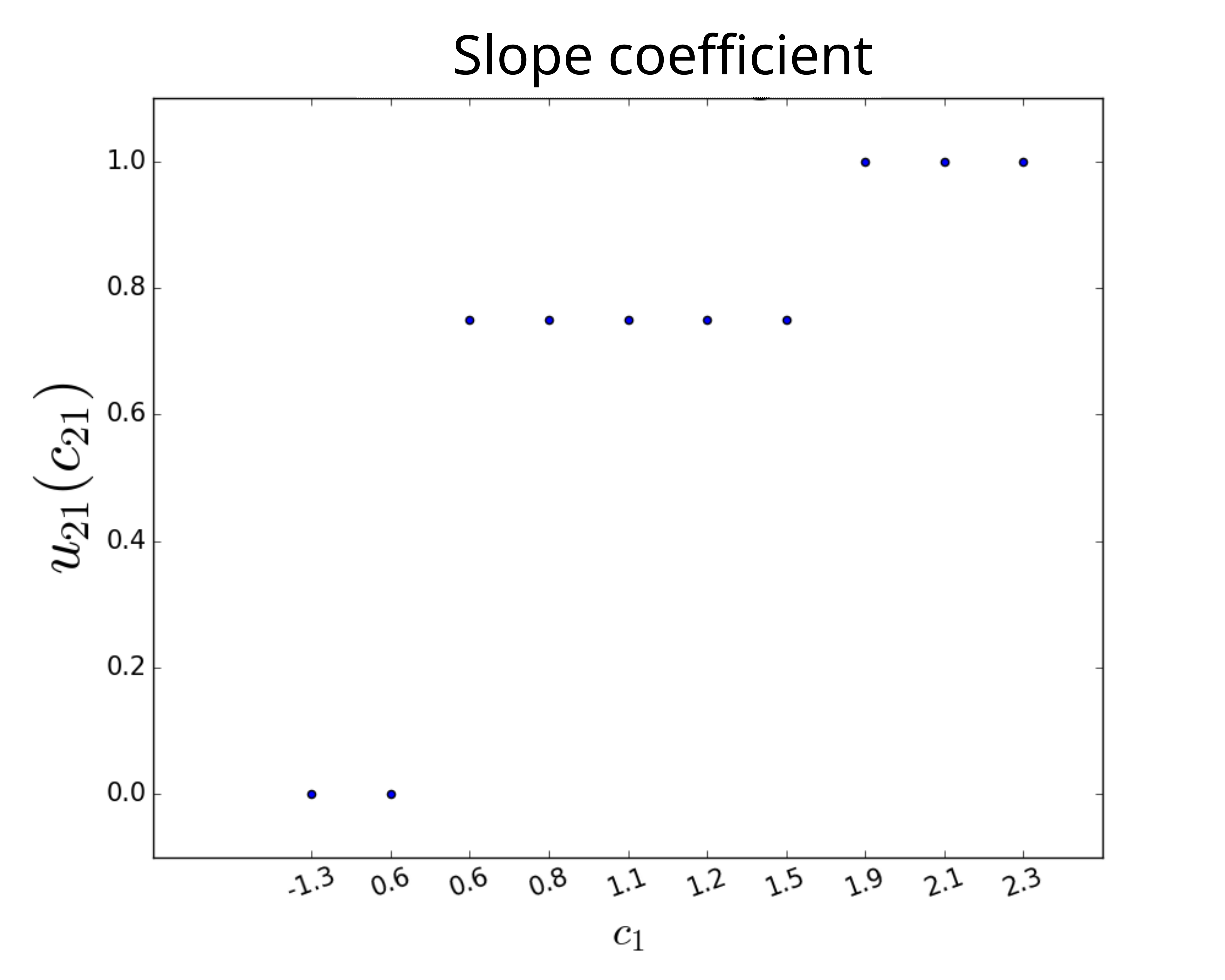}
	\includegraphics[width=6cm]{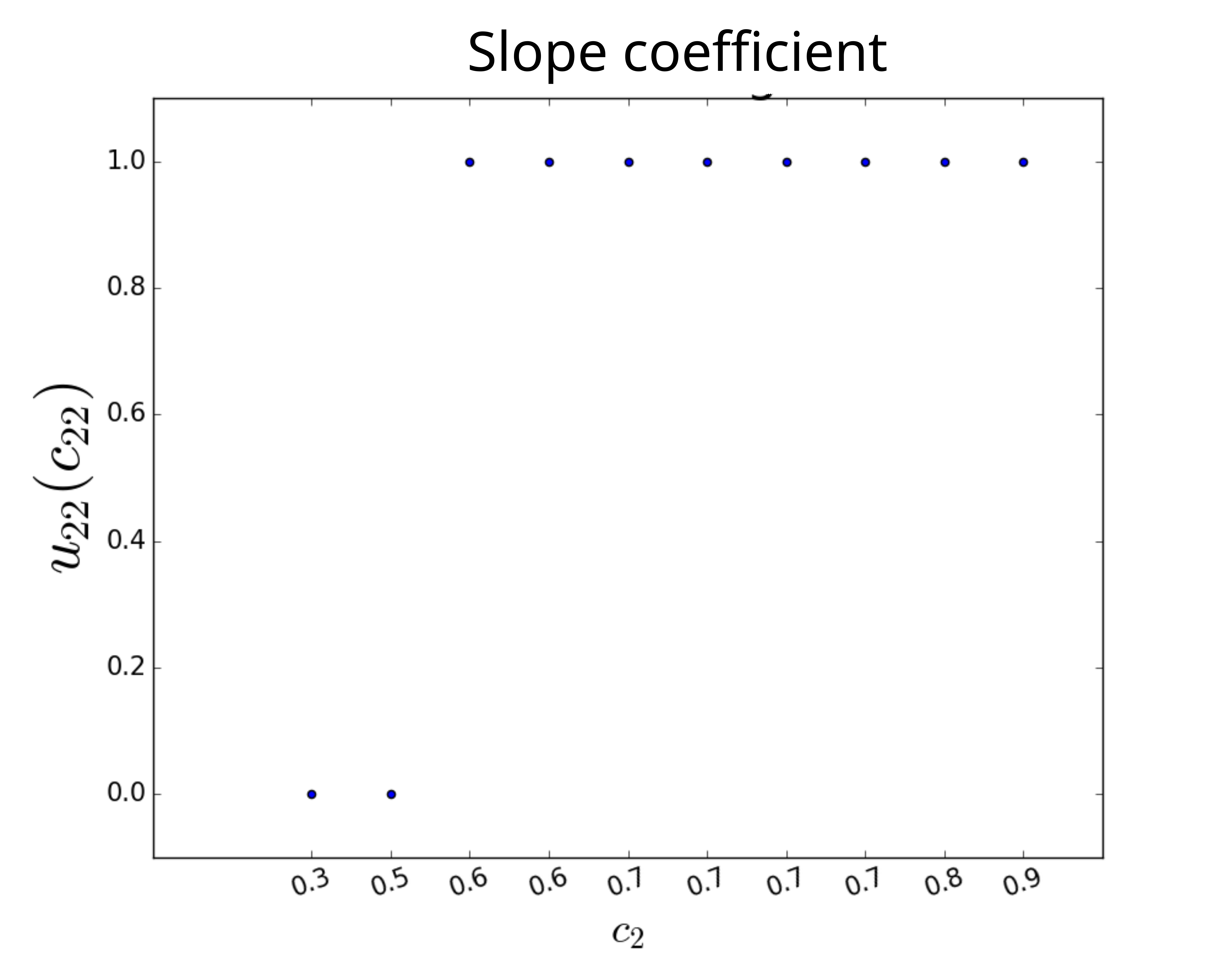}
	\includegraphics[width=6cm]{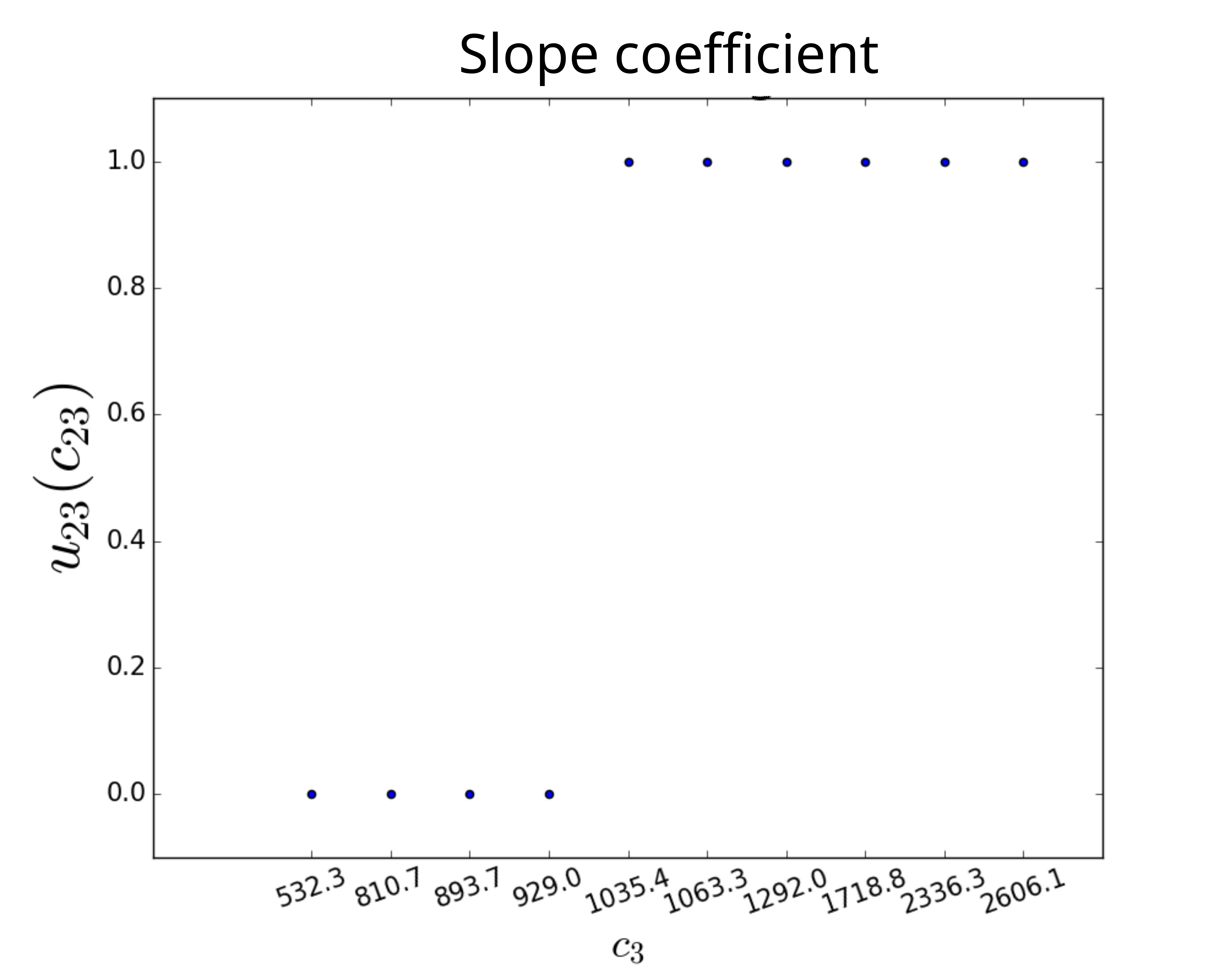}
	\caption{Normalized value functions of descriptive measure 2: slope coefficient.}\label{fig:k2utilitygraph}
\end{figure}

Our goal is to analyze the decision-maker's preferences based on the given ranking, by applying the UTASTAR-T method. Initially, we apply \textbf{Step 1} of UTASTAR-T and obtain descriptive measures, as presented in Table~\ref{tab:dados}. Subsequently, by implementing \textbf{Step 2} to \textbf{Step 4}, an optimal solution is achieved: objective function $z^*$ = 0; and the variables $w_{kj \ell}$ are presented in Table~\ref{tab:variablesresults}. It is important to note that, as $z^* = 0 \implies$ $\sigma^+(a_i) = \sigma^-(a_i) = 0$, which means that there is at least one solution that explains the decision-maker's preference. From the $w_{kj\ell}$ results in Table~\ref{tab:variablesresults}, we obtain the marginal value functions. These marginal values can be normalized by dividing each value function, $u_{kj}(c_{kj}^\ell)$, by the best value of the value function, $u_{kj}(c_{kj}^*)$. The normalized marginal value functions are presented in Figures~\ref{fig:k1utilitygraph} and~\ref{fig:k2utilitygraph}. In addition, the total value of each alternative is:

\begin{minipage}[t]{0.5\textwidth}
	u[\textbf{c}(MY)] =  1,80 
	
	u[\textbf{c}(RU)] =  1,20 
	
	u[\textbf{c}(TR)] =   1,15 
	
	u[\textbf{c}(BR)] =  1,10 
	
	u[\textbf{c}(CN)] =   0,90

\end{minipage}
\begin{minipage}[t]{0.5\textwidth}
	u[\textbf{c}(IN)] =   0,85
	
	u[\textbf{c}(ID)] =  0,80 
	
	u[\textbf{c}(MX)] =   0,75 
	
	u[\textbf{c}(PH)] =  0,70 
	
	u[\textbf{c}(ZA)] =   0,00
\end{minipage}

\vspace{0.2cm}

From Table~\ref{tab:variablesresults} we can observe that it is possible to analyze the decision-maker's preferences in terms of descriptive measures and criteria. For example, we can see that $c_2$ (education) had less relevance in the ranking than $c_1$ (life expectancy) and $c_3$ (gross national income per capita), since its weight is lower for all descriptive measures. Also, $c_2$ (education) is more relevant in terms of trend than in terms of average. The average of $c_3$ (gross national income per capita) is more relevant than that of $c_1$ (life expectancy), but the trend of $c_1$ (life expectancy) is more relevant than that of $c_3$ (gross national income per capita). The same analysis can be made from the value functions shown in Figures~\ref{fig:k1utilitygraph} and~\ref{fig:k2utilitygraph}.

In addition, Figures~\ref{fig:k1utilitygraph} and~\ref{fig:k2utilitygraph} also allow for analyzing the behavior of the value curves for each descriptive measure and each criterion. For example, in terms of average, $c_1$ initially presents null values for the three lowest values and the maximum value for higher values. For $c_3$, the values are null for almost all values, except for the last three.

\begin{table}[!h]
	\centering
	\caption{Results of the variables $w_{kj\ell}$ and their occurrences obtained in the post-optimization multi-objective simulation analysis.}
	\begin{tabular}{cclllllllll|c}
		\toprule
		$k$ & $c_j$&$w_{kj\ell}$ & \multicolumn{8}{c}{Value} & WA  \\ \hline
		\multicolumn{2}{c}{occurrences:} & &                   97 &                    167 &                     69 &                     156 &                     187 &                     53 &                    109 &                     162 &     = 1.000 \\
		\midrule
		\multirow{11}{*}{Average} &  \multirow{3}{*}{$c_1$}& $w_{111}$ &                    0.0 &                    0.0 &                    0.0 &                     0.0 &    \textcolor{red}{0.1} &  \textcolor{red}{0.55} &   \textcolor{red}{0.7} &    \textcolor{red}{0.7} &  0.24 \\
		&&$w_{112}$ &                    0,0 &                    0,0 &                    0,0 &                     0,0 &                     0,0 &   \textcolor{red}{0,1} &                    0,0 &                     0,0 &  0,01 \\
		&&$w_{113}$ &                    0,0 &                    0,0 &                    0,0 &                     0,0 &    \textcolor{red}{0,1} &                    0,0 &                    0,0 &  \textcolor{red}{0,025} &  0,02 \\\cdashline{2-12}
		&  \multirow{5}{*}{$c_2$}&$w_{123}$ &                    0,0 &                    0,0 &                    0,0 &  \textcolor{red}{0,325} &  \textcolor{red}{0,275} &                    0,0 &                    0,0 &                     0,0 &  0,10 \\
		&&$w_{124}$ &                    0,0 &                    0,0 &                    0,0 &   \textcolor{red}{0,15} &                     0,0 &                    0,0 &                    0,0 &                     0,0 &  0,02 \\
		&&$w_{126}$ &                    0,0 &                    0,0 &                    0,0 &                     0,0 &                     0,0 &                    0,0 &  \textcolor{red}{0,05} &                     0,0 &  0,01 \\
		&&$w_{127}$ &                    0,0 &                    0,0 &                    0,0 &                     0,0 &   \textcolor{red}{0,15} &  \textcolor{red}{0,35} &  \textcolor{red}{0,25} &                     0,0 &  0,07 \\
		&&$w_{129}$ &                    0,0 &                    0,0 &                    0,0 &  \textcolor{red}{0,475} &  \textcolor{red}{0,375} &                    0,0 &                    0,0 &                     0,0 &  0,14 \\\cdashline{2-12}
		&  \multirow{3}{*}{$c_3$}&$w_{135}$ &  \textcolor{red}{0,05} &  \textcolor{red}{0,05} &   \textcolor{red}{0,1} &                     0,0 &                     0,0 &                    0,0 &                    0,0 &                     0,0 &  0,02 \\
		&&$w_{137}$ &  \textcolor{red}{0,75} &  \textcolor{red}{0,75} &  \textcolor{red}{0,65} &   \textcolor{red}{0,05} &                     0,0 &                    0,0 &                    0,0 &   \textcolor{red}{0,05} &  0,26 \\
		&&$w_{138}$ &   \textcolor{red}{0,2} &   \textcolor{red}{0,2} &  \textcolor{red}{0,25} &                     0,0 &                     0,0 &                    0,0 &                    0,0 &  \textcolor{red}{0,225} &  0,11 \\ 
		\midrule
		\multirow{13}{*}{SC} &  \multirow{6}{*}{$c_1$}	&$w_{211}$ &                    0,0 &                    0,0 &                    0,0 &                     0,0 &                     0,0 &   \textcolor{red}{0,7} &   \textcolor{red}{0,7} &   \textcolor{red}{0,75} &  0,23 \\
		&&$w_{213}$ &                    0,0 &                    0,0 &                    0,0 &                     0,0 &                     0,0 &  \textcolor{red}{0,05} &  \textcolor{red}{0,05} &   \textcolor{red}{0,05} &  0,02 \\
		&&$w_{216}$ &   \textcolor{red}{0,1} &   \textcolor{red}{0,1} &  \textcolor{red}{0,05} &                     0,0 &                     0,0 &   \textcolor{red}{0,1} &   \textcolor{red}{0,1} &  \textcolor{red}{0,075} &  0,06 \\
		&&$w_{217}$ &                    0,0 &  \textcolor{red}{0,05} &                    0,0 &                     0,0 &                     0,0 &  \textcolor{red}{0,05} &                    0,0 &   \textcolor{red}{0,05} &  0,02 \\
		&&$w_{218}$ &                    0,0 &                    0,0 &                    0,0 &   \textcolor{red}{0,05} &   \textcolor{red}{0,05} &                    0,0 &                    0,0 &                     0,0 &  0,02 \\
		&&$w_{219}$ &                    0,0 &                    0,0 &                    0,0 &                     0,0 &                     0,0 &  \textcolor{red}{0,05} &   \textcolor{red}{0,1} &                     0,0 &  0,01 \\\cdashline{2-12}
		& \multirow{3}{*}{$c_2$}&$w_{222}$ &                    0,0 &                    0,0 &                    0,0 &  \textcolor{red}{0,275} &  \textcolor{red}{0,525} &                    0,0 &                    0,0 &                     0,0 &  0,14 \\
		&&$w_{224}$ &                    0,0 &                    0,0 &                    0,0 &  \textcolor{red}{0,525} &  \textcolor{red}{0,425} &  \textcolor{red}{0,05} &                    0,0 &                     0,0 &  0,16 \\
		&&$w_{225}$ &  \textcolor{red}{0,05} &                    0,0 &  \textcolor{red}{0,15} &                     0,0 &                     0,0 &                    0,0 &                    0,0 &                     0,0 &  0,02 \\\cdashline{2-12}
		&\multirow{4}{*}{$c_3$}&$w_{231}$ &  \textcolor{red}{0,75} &   \textcolor{red}{0,7} &   \textcolor{red}{0,7} &   \textcolor{red}{0,05} &                     0,0 &                    0,0 &                    0,0 &                     0,0 &  0,25 \\
		&&$w_{234}$ &  \textcolor{red}{0,05} &   \textcolor{red}{0,1} &                    0,0 &                     0,0 &                     0,0 &                    0,0 &  \textcolor{red}{0,05} &  \textcolor{red}{0,075} &  0,04 \\
		&&$w_{236}$ &                    0,0 &                    0,0 &  \textcolor{red}{0,05} &    \textcolor{red}{0,1} &                     0,0 &                    0,0 &                    0,0 &                     0,0 &  0,02 \\
		&&$w_{239}$ &  \textcolor{red}{0,05} &  \textcolor{red}{0,05} &  \textcolor{red}{0,05} &                     0,0 &                     0,0 &                    0,0 &                    0,0 &                     0,0 &  0,02 \\
		\bottomrule
	\end{tabular}\label{tab:resultadoexpIDHmultiobjetivovariassolucoes}
\end{table}

\begin{table}[!h]
	\centering
	\caption{Result of the variables $w_{kj\ell}$ and their occurrences obtained in the simulation of the multi-objective post-optimization analysis when there is an order of relevance of the criteria.}
	\begin{tabular}{clllllllllll}
		\toprule
		 Occur. & \multicolumn{11}{c}{$c_1 \succ c_2 \succ c_3$ } \\\hline
		&  $w_{111}$ &  $w_{112}$ &  $w_{126}$ &  $w_{127}$ &  $w_{211}$ &  $w_{213}$ &  $w_{216}$ &  $w_{217}$ &  $w_{219}$ &  $w_{224}$ &  $w_{234}$\\
		\hline
		338 &          0,55 &           0,1 &          0,00 &          0,35 &           0,7 &          0,05 &           0,1 &          0,05 &          0,05 &          0,05 &          0,00 \\
		662 &          0,70 &           0,0 &          0,05 &          0,25 &           0,7 &          0,05 &           0,1 &          0,00 &          0,10 &          0,00 &          0,05 \\	
		\midrule 
		  Occur. & \multicolumn{11}{c}{$c_1 \succ c_3 \succ c_2$} \\\hline
		&  $w_{111}$ &  $w_{113}$ &  $w_{137}$ &  $w_{138}$ &  $w_{211}$ &  $w_{213}$ &  $w_{216}$ &  $w_{217}$ &  $w_{234}$ & & \\
		\hline
		1000 &           0,7 &         0,025 &          0,05 &         0,225 &          0,75 &          0,05 &         0,075 &          0,05 &         0,075 & & \\
		\midrule 
		 Occur. & \multicolumn{11}{c}{$c_2 \succ c_1 \succ c_3$} \\\hline	
		&  $w_{111}$ &  $w_{113}$ &  $w_{123}$ &  $w_{127}$ &  $w_{129}$ &  $w_{218}$ &  $w_{222}$ &  $w_{224}$& & & \\
		\hline
		1000 &           0,1 &           0,1 &         0,275 &          0,15 &         0,375 &          0,05 &         0,525 &         0,425 && & \\
		\midrule
		  Occur. & \multicolumn{11}{c}{$c_2 \succ c_3 \succ c_1$ } \\\hline
		&  $w_{123}$ &  $w_{124}$ &  $w_{129}$ &  $w_{137}$ &  $w_{218}$ &  $w_{222}$ &  $w_{224}$ &  $w_{231}$ &  $w_{236}$& &  \\
		\hline
		1000 &         0,325 &          0,15 &         0,475 &          0,05 &          0,05 &         0,275 &         0,525 &          0,05 &           0,1& &  \\
		\midrule
		 Occur. & \multicolumn{11}{c}{$c_3 \succ c_1 \succ c_2$ } \\\hline
		&  $w_{135}$ &  $w_{137}$ &  $w_{138}$ &  $w_{216}$ &  $w_{217}$ &  $w_{231}$ &  $w_{234}$ &  $w_{239}$  & & & \\
		\hline
		1000 &          0,05 &          0,75 &           0,2 &           0,1 &          0,05 &           0,7 &           0,1 &          0,05 && &  \\
		\midrule   
		  Occur. & \multicolumn{11}{c}{ $c_3 \succ c_2 \succ c_1$} \\\hline
		&  $w_{135}$ &  $w_{137}$ &  $w_{138}$ &  $w_{216}$ &  $w_{225}$ &  $w_{231}$ &  $w_{234}$ &  $w_{236}$ &  $w_{239}$ & & \\ \hline
		678 &          0,05 &          0,75 &          0,20 &          0,10 &          0,05 &          0,75 &          0,05 &          0,00 &          0,05& &  \\
		322 &          0,10 &          0,65 &          0,25 &          0,05 &          0,15 &          0,70 &          0,00 &          0,05 &          0,05 & & \\
		\bottomrule
	\end{tabular}\label{tab:resultadoexpIDHmultiobjetivocompreferenciadoscriterios}
\end{table}

\begin{figure}[h]
	\centering
	\includegraphics[width=6cm]{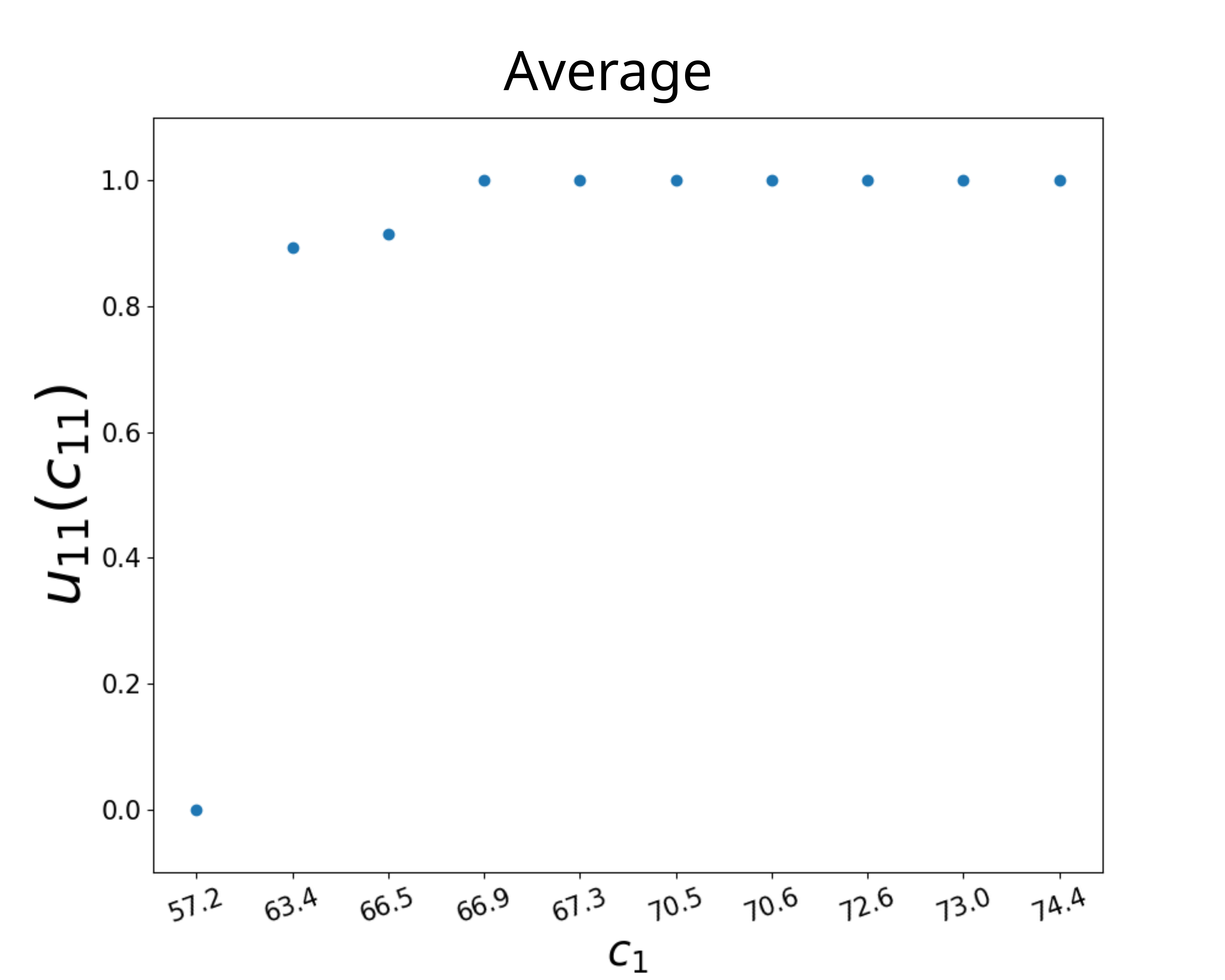}
	\includegraphics[width=6cm]{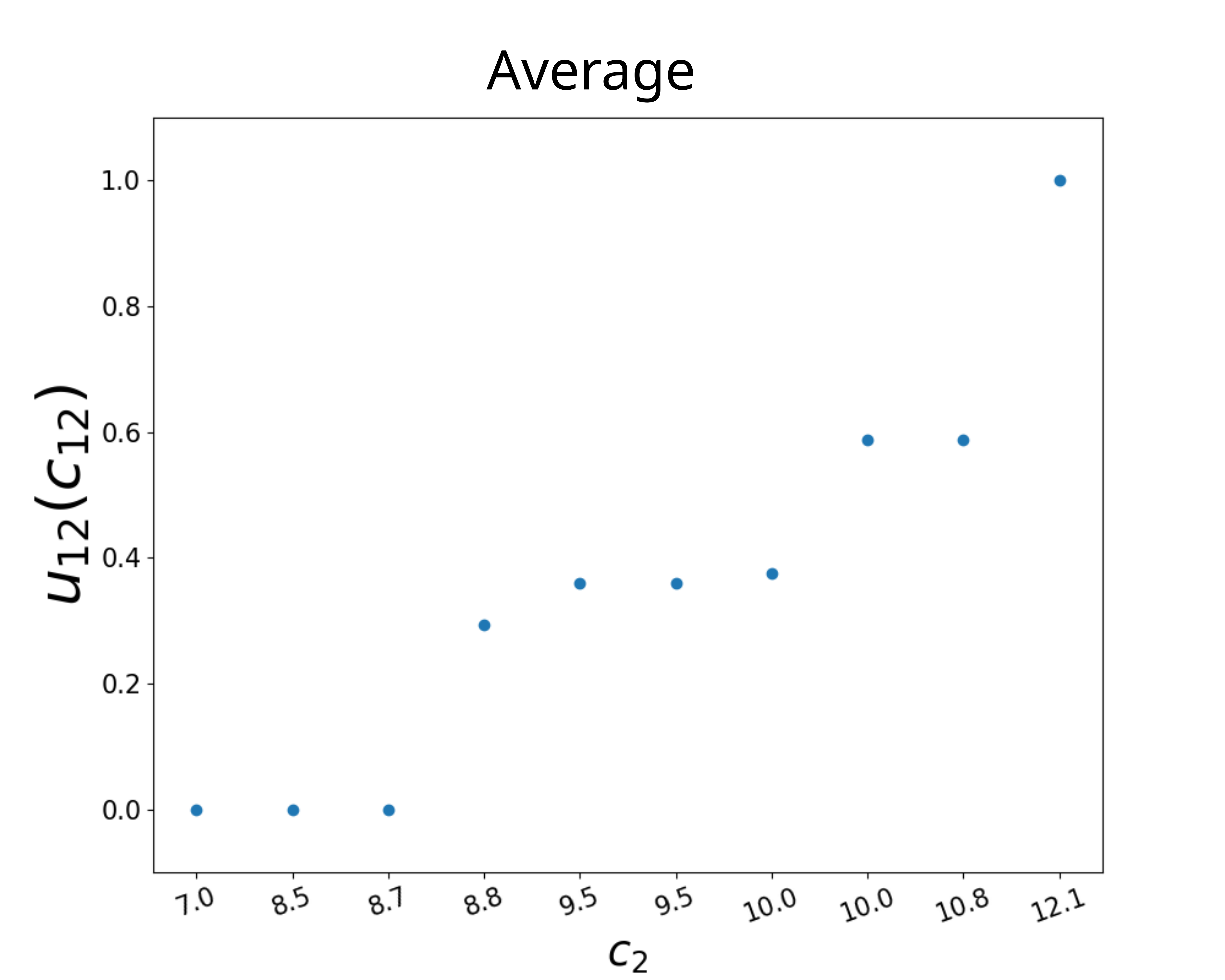}
	\includegraphics[width=6cm]{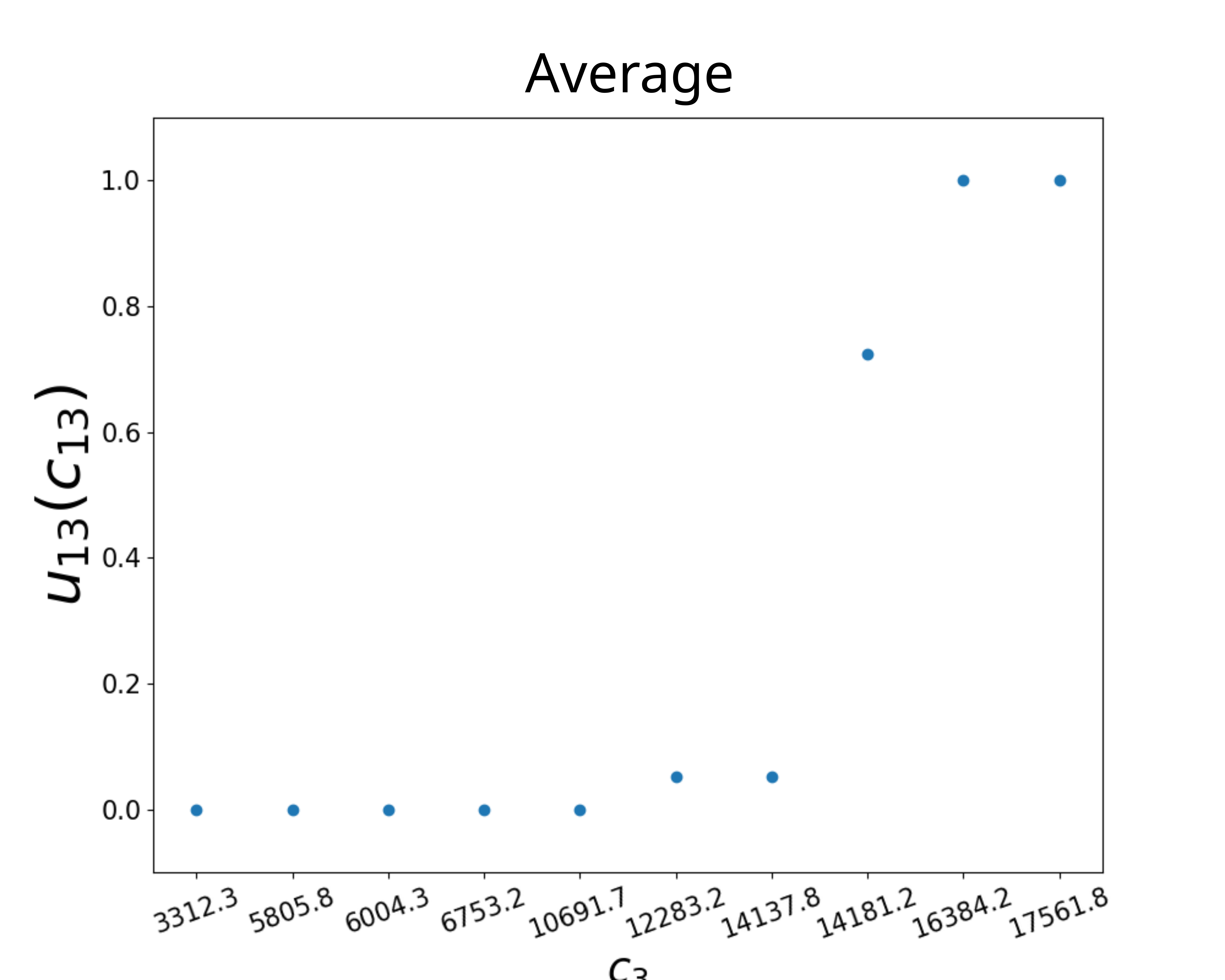}
	\caption{Weighted average of value function values for descriptive measure 1: average.}\label{fig:k1utilitygraph_multi}
\end{figure}
\begin{figure}[h]
	\centering
	\includegraphics[width=6cm]{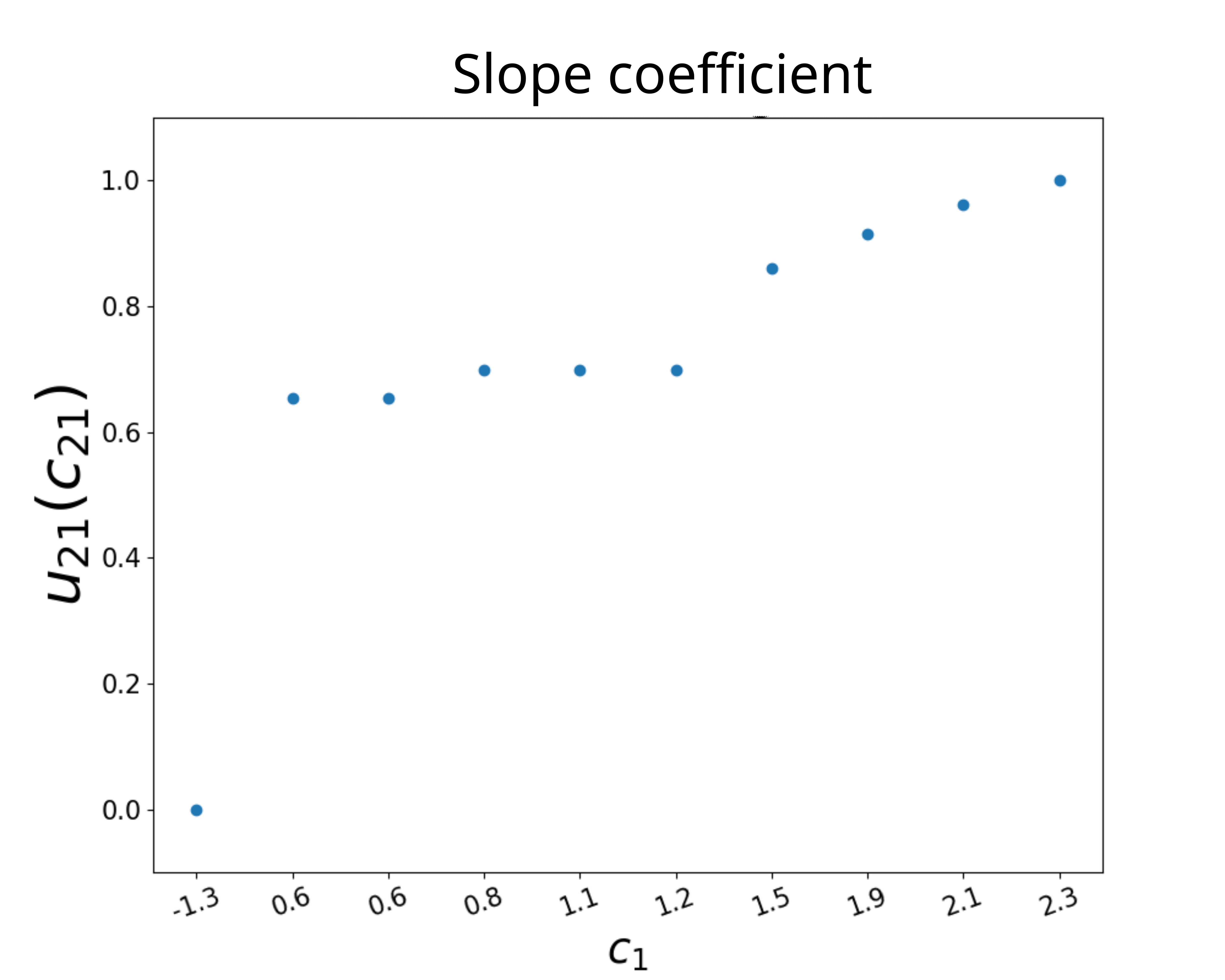}
	\includegraphics[width=6cm]{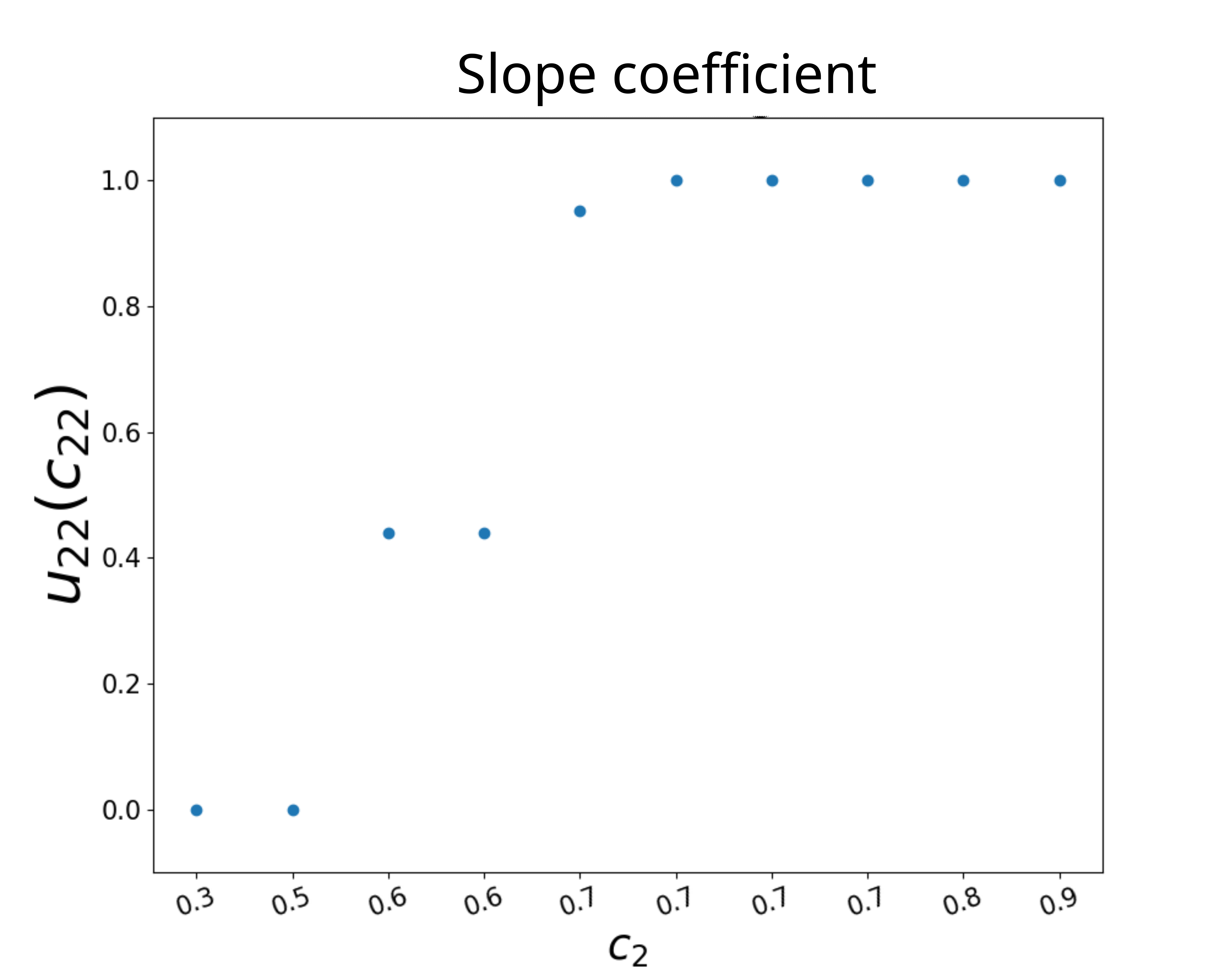}
	\includegraphics[width=6cm]{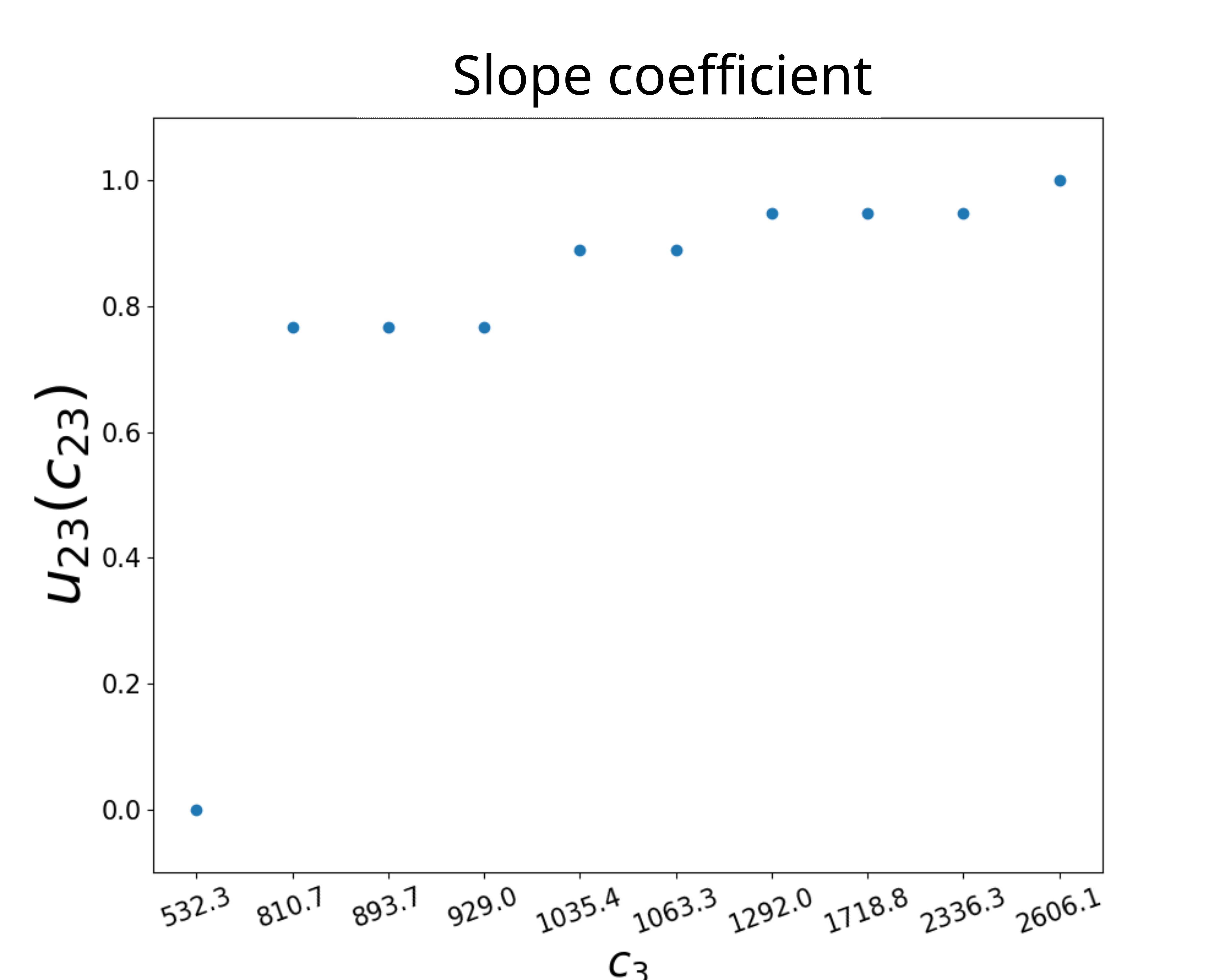}
	\caption{Weighted average of value function values for descriptive measure 2: slope coefficient.}\label{fig:k2utilitygraph_multi}
\end{figure}

For the multi-objective analysis, we developed the Tables~\ref{tab:resultadoexpIDHmultiobjetivovariassolucoes} and~\ref{tab:resultadoexpIDHmultiobjetivocompreferenciadoscriterios} and the figures~\ref{fig:k1utilitygraph_multi} and~\ref{fig:k2utilitygraph_multi}. Table~\ref{tab:resultadoexpIDHmultiobjetivovariassolucoes} shows the results of the multi-objective post-optimization analysis,\textbf{ Step MO-1} to \textbf{Step MO-3}. The first column indicates the characteristic $k$, the second the criterion $c_j$, and the third the values of $w_{kj\ell}$ that were different from zero in at least one result of the simulation. The next eight columns show the results of the $w_{kj\ell}$ and how many times that same result occurred in the simulation. Finally, the last column shows the weighted average (WA) of the simulation's $w_{kj}$ values. Figures~\ref{fig:k1utilitygraph_multi} and~\ref{fig:k2utilitygraph_multi} illustrate the results of the weighted average (from Table~\ref{tab:resultadoexpIDHmultiobjetivovariassolucoes}) in the form of graphs. Table~\ref{tab:resultadoexpIDHmultiobjetivocompreferenciadoscriterios} shows the results of the multi-objective post-optimization analysis, now with the \textbf{Step MO-1 - With criteria ranking} and \textbf{Step MO-2} and \textbf{Step MO-3}, for each criteria preference order.

One initial observation documented in Table~\ref{tab:resultadoexpIDHmultiobjetivovariassolucoes} is that the values of each $w_{kj\ell}$ can exhibit significant variation. For example, the variable $w_{111}$ can have its value ranging from 0 to 0.7, or $w_{137}$ which varies between 0.75 and 0. To better understand this variation of $w_{kj\ell}$, one can analyze $w_{111}$ and $w_{137}$ in Table~\ref{tab:resultadoexpIDHmultiobjetivocompreferenciadoscriterios}, we notice that the values of $w_{111}$ are only positive and take on a relatively high value when $c_1$ is the most relevant criterion. When we have $c_2 \succ c_1 \succ c_3$, $w_{111} = 0.10$ and, in all other cases, $w_{111} = 0$. As for the variable $w_{137}$, it is observed to be positive and takes on a relatively high value when $c_3$ is the most relevant criterion while having little to no significance in all other cases. Thus, this large variation of $w_{kj\ell}$ in Table~\ref{tab:resultadoexpIDHmultiobjetivovariassolucoes} is due to the order in which the $\mu_j$ are obtained in the multi-objective proposal's: \textbf{Step MO-1}.

Therefore, we realize that the values of the $w_{kj\ell}$ in Table~\ref{tab:resultadoexpIDHmultiobjetivovariassolucoes} are somewhat "grouped" as in Table~\ref{tab:resultadoexpIDHmultiobjetivocompreferenciadoscriterios}. Note from Table~\ref{tab:resultadoexpIDHmultiobjetivovariassolucoes} that the $w_{kj\ell}$ of $c_1$ assume positive and relatively high values in the last 3 columns, with occurrences 53, 109 and 162, totaling 324 occurrences (approximately 1/3 of the total occurrences); the $w_{kj\ell}$ of $c_2$ have positive and relatively high values in the columns with occurrences 156 and 187, totaling 343 occurrences (approximately 1/3 of the total occurrences); the $w_{kj\ell}$ of $c_3$ have positive and relatively high values in the columns with occurrences 97, 167, 69, totaling 333 occurrences (approximately 1/3 of the total occurrences). From this analysis, we can conclude that, for the data analyzed, the values of the $w_{kj\ell}$ are very sensitive to which value function is being maximized, which shows that obtaining the upper and lower limits and taking the average may not be representative depending on the problem. In addition we see that there is a great loss of information if we only consider the weighted average.

The analysis of Table~\ref{tab:resultadoexpIDHmultiobjetivocompreferenciadoscriterios} shows that, if it is possible to understand the decision-maker's preferences regarding the relevance of the criteria, then it is possible to find convergence in the results. In some cases there is 100\% convergence of the values of $w_{kj\ell}$, such as when $c_1 \succ c_3 \succ c_2$, because in the 1,000 simulations, the values of $w_{kj\ell}$ were always the same. It is interesting to note that this convergence only occurs with one order of the $\mu_j$, regardless of the values they can assume. In addition, even in cases where there was no 100\% convergence, such as in $c_1 \succ c_2 \succ c_3$, the values are very similar, there is no considerable difference between them.

In the end, we observed that Figures~\ref{fig:k1utilitygraph_multi} and~\ref{fig:k2utilitygraph_multi}, which illustrate the weighted average of the results, are very different from those obtained only in the UTASTAR-T optimization, Figures~\ref{fig:k1utilitygraph} and~\ref{fig:k2utilitygraph} (without considering the multi-objective model). Thus, we can conclude that considering only the solution of the UTASTAR-T model, without post-optimization analysis, is not representative of the problem.

\section{Conclusions}
\label{sec:conclusions}

In this contribution, we proposed to analyze preferences in decisions involving multiple criteria from a tensorial approach. The main novelty is to obtain decision-maker's value functions considering descriptive measures of time series associated with the evolution of criteria. We proposed the UTASTAR-T method for tensor disaggregation and a multi-objective method to address issues related to post-optimization analysis.

Our computational results demonstrate that decision-maker's preferences can be analyzed in terms of both descriptive measures and criteria. Specifically, we observed that for a given criterion, decision-maker's preferences may differ in terms of both its mean value and trend.

Regarding the proposed multi-objective analysis for post-optimization, we found that it is possible to obtain different solution vectors, sometimes with significantly different values. We concluded that, depending on the problem, taking the average of the upper and lower bounds of value functions is not representative of decision-maker's preferences. This highlights the importance of multi-objective analysis as it allows for obtaining multiple solutions that enable discussion between the analyst and the decision-maker to determine the one that best fits their preferences. Furthermore, the analysis shows that if the decision-maker provides a preference order among criteria, better convergence of solutions can be achieved.

Finally, we believe that the proposed disaggregation regarding descriptive measures allows for a more comprehensive analysis of decision-maker's preferences. This analysis can be extended to include more descriptive measures, if necessary.  Furthermore, since only a few MCDA methods have been developed to effectively analyze data with multiple dimensions, our proposed approach represents a potential avenue for future research and practical applications in real-world decision-making scenarios.

\bibliography{bibliografia}

\end{document}